\documentclass[article]{elsarticle}

\usepackage{
hyperref}

\journal{ESAIM COCV
}

\usepackage{amsfonts}
\usepackage{amsmath}
\usepackage{amssymb}

\usepackage{graphicx}
\usepackage{psfrag}
\usepackage{color}
\usepackage{caption}
\usepackage{mathrsfs}
\usepackage[sans]{dsfont}
\usepackage{hyperref}  					 

\newtheorem{thm}{Theorem}
\newtheorem{lem}{Lemma}[section]
\newtheorem{cor}{Corollary}[section]
\newtheorem{prop}{Proposition}[section]
\newtheorem{defn}{Definition}[section]
\newtheorem{rem}{Remark}[section]
\newtheorem{ex}{Example}[section]

\newcommand{\R}{\mathbb{R}}
\newcommand{\N}{\mathbb{N}}

\newcommand{\ve}{\varepsilon}
\newcommand{\n}{\noindent}

\newcommand{\vt}{\vartheta}
\newcommand{\s}{\sigma}

\newcommand{\mE}{\mathcal{E}}

\newcommand{\mO}{\mathcal{O}}

\begin{document}

\begin{frontmatter}

\title{
Multiplicative controllability for nonlinear degenerate parabolic equations between sign-changing states
\tnoteref{mytitlenote}
}
\tnotetext[mytitlenote]{This work was supported by the Istituto Nazionale di Alta Matematica (INdAM), through the GNAMPA Research Project 2016 
\lq\lq 
Controllo, regolarit\`a e viabilit\`a per alcuni tipi di equazioni diffusive'' 
(coordinator P. Cannarsa), 
and the GNAMPA Research Project 2017 
\lq\lq 
Comportamento asintotico e controllo di equazioni di evoluzione non lineari'' 
(coordinator C. Pignotti). 
Moreover, this research was 
performed in the framework of the GDRE CONEDP (European Research Group on \lq\lq Control of Partial Differential Equations'') issued by CNRS, INdAM and Universit\'e de Provence.
This paper was also supported by the research project of the University of Naples Federico II: \lq\lq Spectral and Geometrical Inequalities''.
}

\author[mymainaddress]{G. Floridia \corref{mycorrespondingauthor}}
\cortext[mycorrespondingauthor]{Corresponding author}
\ead[url]{http://wpage.unina.it/giuseppe.floridia}
\ead{giuseppe.floridia@unina.it \& floridia.giuseppe@icloud.com}

\author[mymainaddress]{C. Nitsch}
\ead[url]{http://wpage.unina.it/c.nitsch}
\ead{c.nitsch@unina.it}

\author[mymainaddress]{C. Trombetti}
\ead[url]{http://wpage.unina.it/cristina/Homepage.html}
\ead{cristina.trombetti@unina.it}

\address[mymainaddress]{Department of Mathematics and Applications \lq\lq R. Caccioppoli'',\\ 
      University of Naples Federico II,
80126 Naples, Italy}


\begin{abstract}
In this paper we study the global approximate multiplicative controllability for nonlinear degenerate parabolic Cauchy problems. In particular, we consider a one-dimensional semilinear degenerate reaction-diffusion equation in divergence form governed via the coefficient of the \-reaction term (bilinear or multiplicative control). The above one-dimensional equation is degenerate since the diffusion coefficient is positive on the interior of the spatial domain and vanishes at the boundary points. Furthermore, two different kinds of degenerate diffusion coefficient are distinguished and studied in this paper: the weakly degenerate case, that is, if the reciprocal of the diffusion coefficient is summable, and the strongly degenerate case, that is, if that reciprocal isn't summable. In our main result we show that the above systems can be steered from an initial continuous state that admits a finite number of points of sign change to a target state with the same number of changes of sign in the same order. Our method uses a recent technique introduced for uniformly parabolic equations employing the shifting of the points of sign change by making use of a finite sequence of initial-value pure diffusion pro\-blems. Our interest in degenerate reaction-diffusion equations is motivated by the study of some \-energy balance models in climatology (see, e.g., the Budyko-Sellers model) and some models in population genetics (see, e.g., the Fleming-Viot model).
%
\end{abstract}


\end{frontmatter}


\section{Introduction}
This paper is concerned with the study 
 of the global approximate controllability of one-dimensional semilinear degenerate reaction-diffusion equations
governed in the bounded domain $(-1,1)$ by means of the \textit{bilinear controls} $\alpha (x,t),$ of the form
\begin{equation}
\label{Psemilineare}
\left\{\begin{array}{l}
\displaystyle{u_t-(a(x) u_x)_x =\alpha(x,t)u+ f(x,t,u)\,\quad \mbox{ in } \; Q_T \,:=\,(-1,1)\times(0,T)}\\ [2.5ex]
\displaystyle{
\begin{cases}
\begin{cases}
\beta_0 u(-1,t)+\beta_1 a(-1)u_x(-1,t)= 0 \quad \;\:t\in (0,T)\, \\
\qquad\qquad\qquad\qquad\qquad\qquad\qquad\qquad\qquad\qquad\quad\quad(\mbox{for }\, WDP)\\
\gamma_0\, u(1,t)\,+\,\gamma_1\, a(1)\,u_x(1,t)= 0
 \qquad\quad\, t\in (0,T)\,
\end{cases}
\\
\quad a(x)u_x(x,t)|_{x=\pm 1} = 0\,\,\qquad\qquad\qquad\;\;\, t\in(0,T)\;\;\quad(\mbox{for }\, SDP)
\end{cases}
}\\ [2.5ex]
\displaystyle{u(x,0)=u_0 (x) \,\qquad\qquad\qquad\qquad\quad\qquad\qquad\;\; \,x\in(-1,1)}~.
\end{array}\right.
\end{equation}
In the semilinear 
Cauchy 
problem
$(\ref{Psemilineare})$
the diffusion coefficient $
a\in C([-1,1])\cap C^1(-1,1),$ positive on $(-1,1),$ 
vanishes at the boundary points of $[-1,1],$ leading to a degenerate parabolic equation.
Furthermore, two kinds of degenerate diffusion coefficient can be distinguished and are studied in this paper. $(\ref{Psemilineare})$ is a weakly degenerate problem $(WDP)$ (see \cite{CF2} and \cite{F2}
) if the diffusion coefficient is such that
$\frac{1}{a}\in L^1(-1,1),$ while the problem $(\ref{Psemilineare})$ is called a strongly degenerate problem $(SDP)$ (see \cite{CFproceedings1} and  \cite{F1}
) if 
$\frac{1}{a}\not\in L^1(-1,1)$.\\
In this paper, we assume that the reaction coefficient, that is the bilinear control $\alpha(x,t),$ is bounded on $Q_T,$ and the initial datum $u_0(x)$ is continuous on the open interval $(-1,1),$ since $u_0$ belongs to the weighted Sobolev space
$H^1_a(-1,1),$
defined as either
\begin{align}\label{H1a}
\{u\in L^2(-1,1)&|\,u \text{ is absolutely continuous in } [-1,1] \text{ and }\;\sqrt{a}\,u_x\in L^2(-1,1)\}\; \text{ for (WDP)},\nonumber\\
\text{ or }\\
\{u\in L^2(-1,1)&|\,u \text{ is locally absol. 
continuous in } (-1,1) \text{ and }\;\sqrt{a}\,u_x\in L^2(-1,1)\}\; \text{ for (SDP)}.\nonumber
\end{align}
See \cite{ACF} for the main functional properties of this kind of weighted Sobolev spaces, in particular we note that the space $H^1_a(-1,1)$ is embedded in $L^\infty(-1,1)$ only in the weakly degenerate case (see also \cite{CFproceedings1}, \cite{CF2}, \cite{F1} 
and \cite{F2}). After introducing in Section \ref{Pbfor} the problem formulation of \eqref{Psemilineare}, in Section \ref{Wp} we recall the main properties of this weighted Sobolev spaces and the well-posedness of \eqref{Psemilineare} (\footnote{We recall that it is well-known (see, e.g., \cite{ACF}) that in the $(WDP)$ case all functions in the domain of the corresponding differential ope\-rator possess a trace on the boundary, in spite of the fact that the operator degenerates at such points. Thus, in the $(WDP)$ case we can consider in \eqref{Psemilineare} the general Robin type boundary conditions. On the other hand, in the $(SDP)$ 
one is forced to restrict to the Neumann type boundary conditions.}).
In the following 
we introduce
the main 
 motivations for studying degenerate parabolic pro\-blems with the above structure.
\subsection*{
Reaction-diffusion equations and their applications}
\noindent 
It is well-known that reaction-diffusion equations 
can be linked to various applied 
models such as chemical reactions, nuclear chain reactions, social sciences and biomedical models.
More ge\-nerally, reaction-diffusion equations or systems 
describe how the concentration of one or more substances 
changes under the influence of some processes such as local 
reactions, where
substances are transformed into each other, and diffusion which causes substances to spread out 
in space. See, e.g., some recent papers by Enrique Zuazua and coauthors for the control theory of reaction-diffusion equations: \cite{TZZ}, \cite{PTZ} and \cite{PZ} (see also \cite{PT}).

In the contest of degenerate reaction-diffusion equations there are several interesting models.
In particular, we recall some models in population genetics, see, e.g., the Fleming-Viot model (for a comprehensive literature of these applications see Epstein's and Mazzeo's book \cite{EM}), and
 some models arising in mathematical finance, see, e.g., the Black-Scholes equation in the theory of option pricing.
Our interest in degenerate parabolic problem \eqref{Psemilineare} is also motivated by the study of
an ener\-gy balance model in climatology: 
the Budyko-Sellers model. This model was introduced, independently, by Budyko 
(see \cite{B2})
and Sellers (see \cite{S}). The Budyko-Sellers model studies
 the role played by
continental and oceanic areas of ice on the evolution of the climate.
A complete treatment of the mathematical formulation of the Budyko-Sellers model has been obtained by I.J. Diaz in \cite{D2} (see also 
\cite{CFproceedings1} and \cite{CMVams16}).

\subsection*{Structure of the paper}
\noindent 
%
In Section \ref{mmsac} we show the mathematical motivations, the state of the art and the contents of this paper. Then,
in Section \ref{Pbfor} we give the problem formulation and in Section \ref{Wp} we recall the well-posedness of \eqref{Psemilineare} (obtained in \cite{F2} for $(WDP)$, and in \cite{F1} for $(SDP)$). Thus, in Section \ref{main} we introduce the main result for system \eqref{Psemilineare}, that is, Theorem~\ref{th1.1}, together with some of its consequences.
In Sections \ref{Sec Control Strategy}--\ref{proof main result}, we explain the iterative structure of the proof of the main result, and, after introducing two necessary technical tools, Theorem~\ref{th2.2} and Theorem~\ref{th preserving}, we proceed with the proof of Theorem~\ref{th1.1}.
Section \ref{ADCP}  deals with the proof of Theorem~\ref{th2.2}: a controllability result for pure diffusion problems. 
Section \ref{nonneg gen} is devoted to the proof of Theorem~\ref{th preserving}: a smoothing result intended to attain suitable intermediate data while preserving the 
points of sign change. 

\subsection{Mathematical motivations, state of the art and contents 
}\label{mmsac}
The interest in degenerate parabolic equations is motivated by several mathematical mo\-dels (see 
Epstein's and Mazzeo's book \cite{EM}
) and is dated back by many decades, in particular significant contributions are due to Fichera's and Oleinik's recherches (see e.g., respectively, \cite{Fichera} and \cite{OR}). In control theory only in the last fifteen years se\-veral contributions about degenerate PDEs appeared, in particular we recall some papers due to Cannarsa and collaborators, see, e.g.,
\cite{ACF}, \cite{ACL} 
and \cite{CMV2}-\cite{CTY
} (principally, we call to mind the pio\-neering and fundamental paper \cite{CMV3}, obtained in collaboration with Martinez and Vancostenoble). 
In the above papers, in \cite{CFF}, \cite{MRV}, and also in many works about controllability for non-degenerate 
equations (see, e.g., 
\cite{FPZ}, \cite{FCZ}, 
 \cite{BFH} and \cite{BDDM1}
), 
boundary and interior locally distributed controls are usually employed, these controls are additive terms in the equation and have localized support. Additive control problems for the Budyko-Sellers model have been studied by J.I. Diaz in \cite{D1} (see also \cite{D2}
).\\
However, such controls are unfit to study several interesting applied problems such as chemical reactions controlled by catalysts, and also smart materials, which are able to change their principal parameters under certain conditions. 
In the present work, the control action takes the form of a bilinear control, that is, a control given by the multiplicative coefficient $\alpha$ in \eqref{Psemilineare}.
General references in the area of 
multiplicative controllability 
are the seminal work \cite{BS} 
by Ball, Marsden, and Slemrod, 
 some important results about bilinear control of the Schr\"odinger equation obtained by Beauchard, Coron, Gagnon, Laurent and Morancey in \cite{BL}, \cite{CGM}, and in the references therein, and some results obtained by Khapalov for parabolic and hyperbolic equations, see \cite{KB}, and the references therein. See also some results for reaction-diffusion equations (both degenerate and uniformly parabolic) obtained by Cannarsa and Floridia in \cite{CFproceedings1} and \cite{CF2}, by Floridia in \cite{F1}, 
 and by Cannarsa, Floridia and Khapalov in \cite{CFK}. Moreover, we mention the recent papers \cite{OTB} and \cite{EB} about multiplicative controllability of heat and wave equation, respectively. New perspectives in the area of the multiplicative controllability are suggested by Enrique Zuazua and coauthors in \cite{TZZ}, in that paper the bilinear control for reaction-diffusion equations has a big applied importance, it represents the so-called Allee threshold.\\
\noindent {\bf Additive vs multiplicative controllability.} Historically, the concept of controllability emerged  in the second
half of the twentieth century in the context of linear ordinary differential equations and was motivated by several engineering, economics and Life sciences applications.
Then, it  was extended to various linear partial differential equations governed by {\it additive} locally distributed (i.e., supported on a bounded subdomain of the space 
domain), lumped (acting at a point), and boundary controls (see, e.g. 
Fattorini and Russell in \cite{FattRuss}, and many papers by J.L. Lions and collaborators
). 
Methodologically, these studies are typically based on the {\em linear} duality pairing technique between  the  control-to-state mapping at hand and its dual observation map (see in \cite{BFH} the Hilbert Uniqueness Method, HUM, introduced by J.L. Lions in 1988), using in some cases Carleman estimates tool (see \cite{ACF}, \cite{CFGY} and \cite{CFY1}).
When this mapping is nonlinear, as it happens in the case of the multiplicative controllability, 
the aforementioned  approach does not apply and the above-stated concept of controllability becomes, in general, unachievable. \\
In the last years, in spite of the mentioned difficulties, 
many researchers started to study multiplicative controllabi\-lity, since additive controls (see also \cite{ACF} and \cite{CMV2}--\cite{CMV17}) are unable to treat application problems that require inputs with high energy levels or they are not available due to the physical nature of the process at hand. Thus, an approach based on multiplicative controls, where the coefficient $\alpha$ in \eqref{Psemilineare} is used to change the main physical cha\-racteristics of the system at hand, seems realistic.\\
\noindent {\bf State of the art for uniformly parabolic equations.} To motivate the multiplicative controllability results obtained in this paper for degenerate equation, we start to present the state of the art for uniformly parabolic equations. Let us introduce the 
 following semilinear Dirichlet boundary value problem, studied in \cite{CFK},
\begin{equation}\label{1.1}
   \begin{cases}
\quad   u_t \; = \; u_{xx} \; + \; v(x,t)  u \; + \; f(u)
&\quad  {\rm in} \;\;\;(0,1) \times (0, T)\,, \;\; 
\;\; T>0,
\\
\quad u (0,t) = u (1,t) = 0,
&\quad\qquad\quad\;\; t \in (0, T),
\\
\quad u\:\mid_{t = 0} \; = u_0,
\end{cases}
\end{equation}
where  $u_0\in H^1_0 ( 0,1 )$ (\footnote{
$H^1_0(-1,1)=\{u\in L^2(-1,1)| u_x\in L^2(-1,1) \text{ and } u(\pm1)=0\}.$
}),  $v \in L^\infty (Q_T) $ is a  bilinear control,
the nonlinear term $f:\R\rightarrow\R$ is assumed to be a Lipschitz function, differentiable at $u=0,$ and satisfying $f(0)=0$.\\
There are some important obstructions to the multiplicative controllability (\footnote{ Let us recall that, in  general terms, an evolution system  is called globally approximately
controllable in a given space $H$ at time
$ T > 0$, if any initial state in $ H$ can be steered  
into any neighborhood of any desirable target state  at time $ T$,
by selecting a suitable control. 
}) of \eqref{1.1}.  We note that system \eqref{1.1} cannot be steered anywhere from the origin.
Moreover, if $u_0 (x) \geq 0$ in $(0,1)$,  then the strong maximum principle (\footnote{ In Remark 2.1 of \cite{CFK}, it was observed that the strong maximum principle for linear parabolic PDEs 
can be extend to the semilinear parabolic system \eqref{1.1}. 
}) demands that the respective solution to \eqref{1.1} remains nonnegative  at any moment of time, regardless of the choice of the bilinear control $ v$. This means that system \eqref{1.1} cannot be steered from any such $ u_0 $ to any target state which is ne\-gative on a nonzero measure set in the space domain.
Owing to the previous obstruction to the multiplicative controllability two kinds of controllabi\-lity are worth studying: nonnegative controllability and controllability between sign-changing states.\\
First, in \cite{KB} Khapalov studied global nonnegative approximate controllability of the one-dimensional 
 nondegenerate semilinear convection-diffusion-reaction equation governed in a bounded domain via bilinear controls. 
Finally, in  \cite{CFK} Cannarsa, Floridia and Khapalov esta\-blished an approximate controllability property for the semilinear system \eqref{1.1}
in suitable classes of functions that change sign, 
not arbitrarily but respecting the structure imposed by the strong maximum principle, like in the seminal paper by Matano~\cite{Matano}.\\
\noindent {\bf State of the art for degenerate parabolic equations: nonnegative controllability.} 
With regard to the degenerate reaction-diffusion equations, similar results about global nonne\-gative approximate multiplicative controllability 
were obtained 
in \cite{CFproceedings1}, \cite{CF2} and \cite{F1} 
(\footnote{ In \cite{CFproceedings1}, \cite{CF2} and \cite{F1}, 
first, the authors proved  that, also in the degenerate case, 
if $u_0\geq 0$ 
then
the respective solution 
remains nonnegative  at any moment of time, regardless of the choice of the bilinear control.}).\\ 
%
At first, Cannarsa and Floridia considered the linear degenerate problem associated to $(\ref{Psemilineare})$  (i.e. when $f\equiv 0$)  
in the two distinct kinds of set-up. 
Namely, in \cite{CF2}
 the \textit{weakly degenerate} linear problem $(WDP)$ (that is, when $\frac{1}{a}\in L^1(-1,1)$) was investigated, and
 in \cite{CFproceedings1} 
 the \textit{strongly degenerate} linear problem $(SDP)$ (that is, when $\frac{1}{a}\not\in L^1(-1,1)$) was studied.
Then, in \cite{F1} Floridia focused on the semilinear strongly degenerate case. 
So, in the three above intermediate steps, studied in the papers \cite{CFproceedings1}, \cite{CF2} and \cite{F1}, the authors obtained global nonnegative ap\-proximate controllability of \eqref{Psemilineare}, 
 {\it in large time}, via bilinear piecewise static controls with initial state $u_0\in L^2(-1,1).$ 
That is,  it has been showed that the above system can be steered {\it in large time}, in
the space of square-summable functions, from any nonzero, nonnegative initial state
into any neighborhood of any desirable nonnegative target-state by bilinear
piecewise  static controls.

\noindent{\bf Contents of the paper: controllability of \eqref{Psemilineare} between sign-changing states.} In this paper, we study 
the multiplicative controllability of the semilinear degenerate reaction-diffusion system \eqref{Psemilineare} when both the initial and target states admit a finite number of points of sign change, in particular we extend to the degenerate settings the results obtained in \cite{CFK} for uniformly parabolic equations.
There are some substantial differences with respect to the work \cite{CFK}.  
 The main technical difficulty to overcome with respect to the uniformly parabolic case, 
is the fact that functions in $H^1_a(-1,1)$ need not be necessarily bounded when the operator is strongly degenerate.
\noindent The $(WDP)$ case is somewhat similar to the uniformly parabolic case,
however in our control problem we have to cope with further difficulties given by the general
\textit{Robin} boundary conditions.
%
 In spite of the above difficulties we are able to prove for the degenerate problem \eqref{Psemilineare} 
 that given an initial datum $ u_0 \in H_a^1 (-1,1)$ with a finite number of changes of sign,
 any target state $ u^*\in H_a^1 (-1,1)$, with as many changes of sign {\it in the same order} (in the sense of Definition \ref{sameorder}) as the given $ u_0, 
 $ can be approximately reached in the $ L^2 (-1,1)$-norm at some time $ T>0,$ choosing suitable reaction coefficients (see Theorem \ref{th1.1}).\\
We adapt to the degenerate system \eqref{Psemilineare} a technique introduced in \cite{CFK}, for uniformly parabolic equations, employing the shifting of the points of sign change by making use of a finite sequence of 
initial-value pure diffusion problems. In particular, we proceed
by splitting the time interval $[0,T]$ into 
$2N$ time intervals ($N\in\N$ will be determined after the crucial Proposition \ref{P2})
$$[0,T]=[0,S_1]\cup[S_1,T_1]\cup\cdots\cup[T_{N-1},S_N]\cup[S_{N},T_N]\cup[T_N,T],$$
on which two alternative actions are applied. On the even intervals $[S_k,T_k]$ we choose suitable initial data, $w_k\in H^1_a(-1,1)\cap C^{2+\beta}([a_0^*,b_0^*])$ $\left(\text{with suitable } \,[a_0^*,b_0^*]\subset (-1,1)\right)$, in pure diffusion problems ($\alpha\equiv 0$) to move the points of sign change to their desired location, whereas on the odd intervals $[T_{k-1},S_k]$ we use piecewise static multiplicative controls $\alpha_k$ to attain such $w_k$'s as intermediate final conditions. More precisely, on $ 
\bigcup_{k=1}^N[S_k,T_k]$ 
we make use of the boundary problems 
\begin{equation*}
 \left\{\begin{array}{l}
\displaystyle{
\quad   w_t \; = \; (a(x) w_x)_x
\;  + \; f(x,t,w), \;\;\;\;\qquad
 {\rm in} \;\;\;
 \displaystyle (-1,1)\times {\small{\small \bigcup_{k=1}^N}}[S_k,T_k]
}\\ [2.5ex]
\displaystyle{
\begin{cases}
\begin{cases}
\beta_0 w(-1,t)+\beta_1 a(-1)w_x(-1,t)= 0 \quad \;\:t\in 
\displaystyle {\small{\small \bigcup_{k=1}^N}}(S_k,T_k)\, \\
\qquad\qquad\qquad\qquad\qquad\qquad\qquad\qquad\qquad\qquad\qquad\qquad\quad\quad(\mbox{for }\, WDP)\\
\gamma_0\, w(1,t)\,+\,\gamma_1\, a(1)\,w_x(1,t)= 0
 \qquad\quad\, t\in \displaystyle {\small{\small \bigcup_{k=1}^N}}(S_k,T_k)
\end{cases}
\\
\quad a(x)w_x(x,t)|_{x=\pm 1} = 0\,\,\qquad\qquad\qquad\;\;\,\, t\in  \displaystyle {\small{\small \bigcup_{k=1}^N}}(S_k,T_k)
\;\;\quad\;\;\;\;\,(\mbox{for }\, SDP)
\end{cases}
}\\ [2.5ex]
\displaystyle{\quad w\mid_{t = S_k} \; = w_{k}(x),
\qquad\qquad\qquad\qquad\quad\;\; \:\,x\in(-1,1),\quad k = 1, \ldots, N}~.
\end{array}\right.
\end{equation*}
where the $ w_{k}$'s are viewed as control parameters to be chosen to generate suitable curves of sign change, which have to be continued along all the $N$ time intervals $[S_k,T_k]$ until each point has reached the desired final position. 
In order to fill the gaps between two successive $[S_k, T_{k}]$'s, on $[T_{k-1},S_k]$ we construct $\alpha_k$ that steers the solution of 
\begin{equation*}
\left\{\begin{array}{l}
\displaystyle{
\quad   u_t \; = \; 
(a(x) u_x)_x \; + \; \alpha_k(x,t)  u \; + \; f(x,t,u)
\;\;\;\quad  {\rm in} \;\,
(-1, 1) \times [T_{k-1}, S_k
],  }\\ [2.5ex]
\displaystyle{
\begin{cases}
\begin{cases}
\beta_0 u(-1,t)+\beta_1 a(-1)u_x(-1,t)= 0 \quad \;\:t\in 
(T_{k-1}, S_k
)\, \\
\qquad\qquad\qquad\qquad\qquad\qquad\qquad\qquad\qquad\qquad\qquad\qquad\qquad(\mbox{for }\, WDP)\\
\gamma_0\, u(1,t)\,+\,\gamma_1\, a(1)\,u_x(1,t)= 0
 \qquad\quad\, t\in 
 (T_{k-1}, S_k
 )
\end{cases}
\\
\quad a(x)u_x(x,t)|_{x=\pm 1} = 0\,\,\qquad\qquad\qquad\;\;\,\, t\in  
(T_{k-1}, S_k
)
\quad\quad\;\:\quad(\mbox{for }\, SDP)
\end{cases}
}\\ [2.5ex]
\displaystyle{\quad u\mid_{t = T_{k-1}} \; = u_{k-1} + r_{k-1}\in H^1_a(-1,1)
}~,
\end{array}\right.\end{equation*}
 from $u_{k-1}+r_{k-1}$ 
to  $w_k,$
where $ u_{k-1}$ and $ w_{k}$ have  
the same points of sign change, 
 and  $\|r_{k-1}\|_{L^2 (-1, 1)}$ is small. The above result is contained in Theorem 3.\\
The fact that such an iterative process can be completed within a finite number of steps ($2N,$ for suitable $N\in\N$) is an important point of the proof.
Such a point follows from precise estimates which is, essentially, the consequence of the following facts: 
\begin{enumerate}
\item[(a)] the sum of the distances of each branch of the null set of the resulting solution of \eqref{1.1} from its target points of sign change decreases  at a linear-in-time rate for curves which are still far away from their corresponding target points;
\item[(b)] the  error caused by the possible displacement of points already near their targets is ne\-gligible.
\end{enumerate}
\noindent {\bf Some open problems.} 
In the future, we intend to investigate the multiplicative controllability for both degenerate and uniformly parabolic equations in higher space dimensions on domains with specific geometries (see, e.g., Section 6 in \cite{CFK}). 
Moreover, we would like to extend our approach to study the approximate controllability of the general formulation 
of the Budyko-Sellers differential problem on a compact surface 
 without boundary. 
Finally, 
we would like to extend our approach to other nonlinear systems of parabolic type, such as the systems of fluid dynamics (see, e.g., \cite{KCPF}), and the porous medium equation.

\subsection{Problem formulation 
}\label{Pbfor}
In this paper, we consider the degenerate problem (\ref{Psemilineare})
under the following assumptions:
 \begin{enumerate}
  \item[(A.1)] $u_0 \in H^1_a(-1,1);$ (\footnote{ The definition of the weighted Sobolev space $H^1_a(-1,1)$ is given in \eqref{H1a}.}) 
 \item[(A.2)] $\alpha
  \in L^\infty (Q_T);$ 
\item[(A.3)] $f:Q_T\times\R\rightarrow \R$ is such that
\begin{itemize}
\item $(x,t,u)\longmapsto f(x,t,u)$ is a CarathŽ\'eodory function on $Q_T\times\R,$ (\footnote{ We say that $f:Q_T\times\R\rightarrow \R$ is a CarathŽ\'eodory function on $Q_T\times\R$ if the following properties hold:
\begin{itemize}
\item[$\star$] $(x,t)\longmapsto f(x,t,u)$ is measurable, for every $u\in\R;$ 
\item[$\star$] $u\longmapsto f(x,t,u)$ is a continuous function, for a.e. $(x,t)\in Q_T.$
\end{itemize}
})
\item $u\longmapsto f(x,t,u)$ is differentiable at $u=0,$
\item
$t\longmapsto f(x,t,u)$ is locally absolutely continuous for a.e. $x\in (-1,1),$ for every $u\in\R,$\,
\item there exist constants $\gamma_*\geq 0, \vartheta\geq1$ 
and $\nu
\geq0$
 such that, $\mbox{for a.e.} (x,t)\in Q_T, \forall u,v\in \R,$ we have
\begin{equation}\label{Superlinearit}
|f(x,t,u)|\leq\gamma_*\,|u|^\vartheta, 
\end{equation}
\begin{equation}\label{fsigni}
-\nu
\big(1+|u|^{\vartheta-1}+|v|^{\vartheta-1}\big)(u-v)^2\leq \big(f(x,t,u)-f(x,t,v)\big)(u-v)\leq\nu (u-v)^2,
\end{equation}
\begin{equation*}
  f_t(x,t,u)\,u\geq-\nu\, 
  u^2\,;\; 
( \footnote{
This assumption is necessary and
used only for the well-posedness, obtained in the paper \cite{F2} for the $(WDP)$ case, and in \cite{F1} for the $(SDP)$ case. 
})
\end{equation*}
\end{itemize}
   \item[(A.4)] $a \in C([-1,1])\cap C^1(-1,1)$ is such that
    $$a(x)>0, \,\, \forall \, x \in (-1,1),\quad a(-1)=a(1)=0,$$
    moreover, we have the following two alternative assumptions:
     \begin{itemize}
    \item[$({A.4}_{WD})$]
 if
 $\;\;\displaystyle\frac{1}{a}\in L^1(-1,1),$ in \eqref{Psemilineare} let us consider the 
 Robin boundary conditions, where
  \begin{itemize}
\item[$\star$]
   $\beta_0,\beta_1,\gamma_0,\gamma_1\in\R,\;\beta_0^2+\beta_1^2>0, \;\gamma_0^2+\gamma_1^2>0,$ satisfy the sign condition\\
     $$\beta_0\beta_1 \leq 0 \;\;\;\mbox{ and }\;\;\; \gamma_0\gamma_1 \geq 0;$$
\end{itemize}
  \item[$({A.4}_{SD})$] if $\;\;\displaystyle\frac{1}{a}\not\in L^1(-1,1)$
      (\footnote{ We note that if $a \in C^1([-1,1])$ follows $\frac{1}{a}\not\in L^1(-1,1).$
      }) 
   and the function $\xi_a(x):=\displaystyle\int_0^x \frac{1}{a(s)}ds 
    \in L^{2\vartheta-1
    }(-1,1),$ in \eqref{Psemilineare} let us consider the weighted Neumann boundary conditions.
  \end{itemize}

\end{enumerate}
\begin{rem}\rm
The principal part of the operator in $(\ref{Psemilineare})$ coincides with that of the Budyko-Sellers climate model 
for
$a(x)=1-x^2$.\,
In this case $\frac{1}{1-x^2}\not\in L^1(-1,1),$ but $\xi_a(x)=\frac{1}{2}\log\left(\frac{1+x}{1-x}\right)
  \in L^p(-1,1),$ for every $p\geq 1,$ so this is an example of strongly degenerate equation, while an example of weakly degenerate coefficient is $a(x)=\sqrt{1-x^2}.$ 
 \end{rem}

\begin{rem}\label{rem1}
\rm
We note that the inequalities 
 \eqref{fsigni} imply that 
 $$\big|f(x,t,u)-f(x,t,v)\big|\leq
\nu
(1+|u|^{\vartheta-1}+|v|^{\vartheta-1})|u-v|,\;\;\mbox{ for a.e. } (x,t)\in Q_T, \forall u,v\in \R. \,
$$
Moreover, under the further assumption $u\longmapsto f(x,t,u)$ is locally absolutely continuous respect to $u,$
the inequalities 
 \eqref{fsigni} in hypothesis $(A.3)$ is equivalent to the following conditions on the function $f$
 $$-\nu
\big(1+|u|^{\vartheta-1}\big)\leq f_u(x,t,u)\leq\nu
,\; \mbox{ for a.e.  }  (x,t)\in Q_T, \forall u\in \R,$$
thus
 \begin{equation*}
\big|f_u(x,t,u)\big|\leq
\nu
(1+|u|^{\vartheta-1}), \; \mbox{ for a.e. } (x,t)\in Q_T, \forall u
\in \R.\;\,
\end{equation*}
\end{rem}
{\footnotesize
In order to clarify the previous Remark \ref{rem1}, we note that,
since, for a.e. $(x,t)\in Q_T,$ $f(x,t,u)$ is locally absolutely continuous respect to $u,$ we have\\
$$\big(f(x,t,u)-f(x,t,v)\big)(u-v)=(u-v)\int^{u}_{v}f_\xi(x,t,\xi)d\xi\leq(u-v)\int^{u}_{v}\nu\,d\xi\leq\nu(u-v)^2,$$
$$\big|f(x,t,u)-f(x,t,v)\big|\leq\int^{max\{u,v\}}_{min\{u,v\}}|f_\xi(x,t,\xi)|d\xi\leq\nu\int^{max\{u,v\}}_{min\{u,v\}}(1+|\xi|^{\vartheta-1})d\xi\leq\nu(1+|u|^{\vartheta-1}+|v|^{\vartheta-1})|u-v|,$$
for a.e. $(x,t)\in Q_T,$ for every $u,\,v\in\R.$ }

\begin{ex}
An example of function $f$ that satisfies the assumptions $(A.3)$ is the following
$$f(x,t,u)=c(x,t)\min\{|u|^{\vartheta-1},1\}u-|u|^{\vartheta-1}u,$$
where $c$ 
is a Lipschitz continuous function. 
\end{ex}
\begin{rem}
\label{r1.2}\rm
We note that system \eqref{Psemilineare} cannot be steered anywhere from the origin.
Moreover,  in \cite{F1}
it was proved that 
 if $u_0 (x) \geq 0$ in $(-1,1)$ 
the respective solution to \eqref{Psemilineare} remains nonnegative  at any moment of time, regardless of the choice of $\alpha(x,t)$. This means that system \eqref{Psemilineare} cannot be steered from any such $ u_0 $ to any target state which is negative on a nonzero measure set in the space domain. 
\end{rem}
\subsection{Well-posedness 
}\label{Wp}
The well-posedness of the (SDP) problem \eqref{Psemilineare} under the assumptions $(A.1)-
(A.4_{SD})$ is obtained in \cite{F1}, while the well-posedness of the (WDP) problem \eqref{Psemilineare} under the assumptions $(A.1)-
(A.4_{WD})$ is obtained in \cite{F2}. 
In order to deal with the well-posedness of degenerate problem $(\ref{Psemilineare})$, it is necessary to recall the weighted Sobolev space $H^1_a(-1,1),$ already introduced, and to define the space $H^2_a(-1,1)$ (see also \cite{ACF}, \cite{CFproceedings1}, \cite{F1}, \cite{CF2} and \cite{F2}):
$$H^2_a(-1,1):=\{u\in H^1_a(-1,1)| \, au_x \in H^1 (-1,1)\}.$$
$H^1_a(-1,1)$
and $H^2_a(-1,1)$ are Hilbert spaces with the natural scalar products induced,
respectively, by the following norms
$$\|u\|_{1,a}^2:=\|u\|_{L^2(-1,1)}^2+|u|_{1,a}^2 \mbox{ and } \|u\|_{2,a}^2:=\|u\|_{1,a}^2+\|(au_x)_x\|^2_{L^2(-1,1)},$$
where $|u|_{1,a}^2:=\|\sqrt{a}u_x\|_{L^2(-1,1)}^2$ is a seminorm.
\noindent In the following, we will sometimes use $\|\cdot\|$ 
 instead of
$\|\cdot\|_{L^2(-1,1)}.$ 
\\

In \cite{CMV2}, see Proposition 2.1 (see also the Appendix of \cite{F1} and Lemma 2.5 in \cite{CMP}), the following result is proved. 
 \begin{prop}\label{caratH2}
Let $a\in  C^1([-1,1])$ ($(SDP)$ case). 
For every $u\in H^2_a(-1,1)$ we have
$$\lim_{x\rightarrow\pm1}a(x)u_x(x)
=0
\qquad\qquad\text{ and }\qquad\qquad
au\in H^1_0 (-1,1)\,.
$$ 
 \end{prop}
Proposition \ref{caratH2} motivates, in the $(SDP)$ case, the following definition of the operator $(A,D(A))$.\\ 

Given $\alpha\in L^\infty (-1,1),$ let us introduce the operator $(A,D(A))$
defined by
\begin{equation}\label{DA}
\left\{\begin{array}{l}
\displaystyle{
D(A)=\begin{cases}
\left\{u\in H^2_a (-1,1)\bigg\rvert
\begin{cases}
\beta_0 u(-1)+\beta_1 a(-1)u_x(-1)= 0 
\\
\gamma_0\, u(1)\,+\,\gamma_1\, a(1)\,u_x(1)= 0 
\end{cases}
 \right\}\quad\mbox{for }\, (WDP)\\
\;\;
H^2_a (-1,1)\qquad\qquad\qquad\qquad\qquad\qquad\qquad\qquad\qquad\qquad\mbox{for }\, (SDP)\;\;\;
\end{cases}
}\\ [4.8ex]
\displaystyle{\;\;\;\;A\,u=(au_x)_x+\alpha\,u, \,
\,\, \forall \,u \in D(A)}~.
\end{array}\right.
\end{equation}

\subsubsection*{The Banach spaces ${\cal{B}}(Q_T)$ and ${\cal{H}}(Q_T)$}
\noindent Given $T>0,$ let us define the Banach spaces:
$$\mathcal{B}(Q_T):=C([0,T];L^2(-1,1))\cap L^2(0,T;H^1_a (-1,1))$$
with the following norm
\begin{equation*}
\label{normaB}
 \|u\|^2_{\mathcal{B}(Q_T)}= \sup_{t\in
[0,T]}\|u(\cdot,t)\|^2
+2\int^T_{0}\int^1_{-1}a(x)u^2_x
dx\,dt\,,
\end{equation*}
and
$$
{\cal{H}}(Q_T):=L^{2}(0,T;D(A)
)\cap H^{1}(0,T;L^2(-1,1))\cap C([0,T];H^{1}_a(-1,1))
$$
with the following norm 
\begin{equation*}
\label{normaH}
\|u\|^2_{\mathcal{H}(Q_T)}=
 \sup_{
 [0,T]}\left(\|u\|
 ^2+\|\sqrt{a}u_x\|
 ^2\right)+\int_0^T\left(\|u_t\|
 ^2+\|(au_x)_x\|
 ^2\right)\,dt.
 \;\;(\footnote{
 It's well known that this norm is equivalent to the Hilbert norm $$\displaystyle\||u|\|^2_{\mathcal{H}(Q_T)}=
 \int_0^T\left(\|u\|
 ^2+\|\sqrt{a}u_x\|
 ^2+\|u_t\|
 ^2+\|(au_x)_x\|
 ^2\right)\,dt.$$
 })
\end{equation*}


Here we give the definition of \lq\lq\textit{strict solutions} to \eqref{Psemilineare}'' (introduced in \cite{F1} for $(SDP)$ and in \cite{F2} for $(WDP)$), that is the notion of solution with initial state belongs to $H^1_a(-1,1).$ 
\begin{defn}\label{strictsolution}
If 
$u_0\in H^1_a(-1,1),$ u is a \textit{strict solution} 
to \eqref{Psemilineare}, 
if $u\in{\cal{H}}(Q_T)$ and
$$
\left\{\begin{array}{l}
\displaystyle{u_t-(a(x) u_x)_x =\alpha(x,t)u+ f(x,t,u)\,\quad\quad\, \mbox{a.e.\;\; in } \; Q_T \,:=\,(-1,1)\times(0,T)}\\ [2.5ex]
\displaystyle{
\begin{cases}
\begin{cases}
\beta_0 u(-1,t)+\beta_1 a(-1)u_x(-1,t)= 0 \quad\;\; \mbox{a.e.\;\;} \;t\in (0,T)\, \\
\qquad\qquad\qquad\qquad\qquad\qquad\qquad\qquad\qquad\qquad\qquad\qquad(\mbox{for }\, WDP)\\
\gamma_0\, u(1,t)\,+\,\gamma_1\, a(1)\,u_x(1,t)= 0
 \qquad\quad\,\mbox{a.e.\;\;} t\in (0,T)\,
\end{cases}
\\
\quad a(x)u_x(x,t)|_{x=\pm 1} = 0\,\,\qquad\qquad\qquad\;\;\,\mbox{a.e.\;\;}\, t\in(0,T)\;\;\,\quad(\mbox{for }\, SDP)\;\;\;
\end{cases}
}\\ [2.5ex]
\displaystyle{u(x,0)=u_0 (x) \,\qquad\qquad\qquad\qquad\quad\qquad\qquad\;\;\;\; \,x\in(-1,1)}~.
\end{array}\right.
(\footnote{
Since $u\in{\cal{H}}(Q_T)\subseteq L^2(0,T;D(A)
)$ we have $u(\cdot,t)\in D(A)$ for $ \text{ a.e. } t\in(0,T)$, 
so we deduce the associated 
boundary condition. 
})
$$
\end{defn}
We proved
the following result in \cite{F1} (see Appendix B) for $(SDP)$ and in \cite{F2} for $(WDP)$.
\begin{prop}\label{exB}
For all $u_0\in H^1_a(-1,1)$ there exists a unique strict solution $u\in{{\cal{H}}(Q_T)}$ to \eqref{Psemilineare}.
\end{prop}
\begin{rem}\label{rem sol}\rm
In \cite{F1} (for the $(SDP)$) and in \cite{F2} (for the $(WDP)$), for initial data in $L^2(-1,1)$ the notion of \lq\lq {\it strong solutions}'' was defined by approximation 
(\footnote{ The notions of \lq\lq {\it strict/strong solutions}'' are classical in PDEs theory, see, for instance, \cite{BDDM1}, pp. 62-64.}). 
In this paper, we consider states continuous on the open interval $(-1,1),$ then we use initial states in $H^1_a(-1,1),$ so consequently we only consider  the notion of \lq\lq\textit{strict solution}''.
\end{rem}

\section{Main results 
}\label{main}
Our main goal is to show that, given a initial datum $ u_0 \in H_a^1 (-1,1)$ with a finite number of changes of sign,
 any target state $ u^* \in H_a^1 (-1,1)$, with as many changes of sign {\it in the same order} 
 (see Definition \ref{sameorder}) 
 as the given $ u_0, 
 $ can be approximately reached in the $ L^2 (-1,1)$-norm at some time $ T>0,$ choosing suitable reaction coefficients. 
Thanks to this result, in Corollary \ref{th1.2} we easily show, by approximation argument, 
that the system \eqref{Psemilineare} can be also steered toward the target states such that the amount of points of sign change is no more than the one of the given initial data. 
Now, we give some definitions to clarify and simplify the notation.
\begin{defn}
We say that $\bar{u}\in H_a^1 (-1,1)$ 
has $n$ {\it points of sign change},
if there exist $n$ points
  $\bar{x}_l,\;l=1,\dots,n,$ with
 $$-1
 <
 \bar{x}_1<\cdots<\bar{x}_n<
 1
 $$ such that
 \begin{itemize}
 \item[$\star$] $\bar{u}(x)=0,\; x\in(-1,1)\quad\Longleftrightarrow\quad x=\bar{x}_l,\quad l=1,\ldots,n;$
 \item[$\star$] for $l=1,\ldots,n,$ 
  $$\bar{u}(x)\bar{u}(y)<0,\;\;\; \forall x\in \left(\bar{x}_{l-1},\bar{x}_l\right),\, \forall y\in \left(\bar{x}_{l},\bar{x}_{l+1}\right),
   $$
   where let us set $\bar{x}_0:=-1$ and $\bar{x}_{n+1}:=1.$
  \end{itemize}
   \end{defn}
       \begin{defn}\label{sameorder}
We say that $u_0,u^* \in H_{a}^1 (-1,1)$ have the $n$ points of sign change in the same order, if 
 denoting by $x_l^0,\,x_l^*,\, l=1,\ldots,n,$
  the zeros
 of $u_0$ and $u^*,$ respectively, 
 we have
 $$u_0(x)u^*(y)>0,\;\;\; \forall x\in \left(x^0_{l-1},x^0_l\right),\, \forall y\in \left(x^*_{l-1},x_{l}^*\right), \text{ for } l=1,\ldots,n+1,$$
 where let us set $x^0_0=x^*_0=-1$ and $x^0_{n+1}=x^*_{n+1}=1.$
 \end{defn}

\begin{defn}\label{static}
We say that a 
function $\alpha\in L^\infty(Q_T)$ is {\it piecewise static}, if there exist $m\in\N,$ $c_k(x)\in L^\infty(-1,1)$ and $t_k\in [0,T], \,t_{k-1}<t_k,\, k=1,\dots,m$ with $t_0=0 \mbox{ and } t_m=T,$ 
such that $$\alpha(x,t)=c_1(x)\mathds{1}_{[t_{0},t_1]}(t)+\sum_{k=2}^m c_k(x)\mathds{1}_{(t_{k-1},t_k]}(t),$$ where $\mathds{1}_{[t_{0},t_1]}\,  \mbox{  and  }  \,\mathds{1}_{(t_{k-1},t_k]}$ are the indicator function of $[t_{0},t_1]$ and $(t_{k-1},t_k]$, respectively.
\end{defn}

\begin{thm}\label{th1.1} Let $ u_0 \in  H_a^1 (-1,1)$. Assume that $ u_0$ has a finite number of 
points of sign change.
Consider any  $ u^*\in H_a^1 (-1,1)$ which has exactly the same number of
 points of sign change in the same order as $ u_0 $. Then, for any $ \eta > 0$, there exists  $ T=T(\eta,u_0,u^*) > 0$ and  a piecewise static multiplicative control $\alpha=\alpha(\eta,u_0,u^*) \in L^\infty (Q_T) $ such that the respective solution $u$ to $\eqref{Psemilineare}$ satisfies
\begin{equation*}
\| u(\cdot, T) - u^* \|_{L^2 (-1,1)} \; \leq \; \eta.
\end{equation*}
\end{thm}

In Figure $a)$ we explain the statement of Theorem \ref{th1.1}.

\vspace{-1.0cm}
\setlength{\unitlength}{1mm}
\begin{picture}(150,55)(0,0)
\linethickness{1pt}
\put (2,15){\vector(1,0){85}}
\put (10,0){\line(0,1){45}}
\put (40,0){\vector(0,1){45}}
\put (70,0){\line(0,1){45}}
\qbezier(10,20)(16,32)(23,18)
\qbezier(23,18)(34,1)(45,10)
\qbezier(45,10)(65,25)(70,17)
{
\qbezier(10,17)(19,42)(27,19)
\qbezier(27,19)(34,-15)(48,18)
\qbezier(48,18)(65,55)(70,20)
}
\put(23,10){\makebox(0,0)[b]{$x^0_{1}$}}
\put(54,10){\makebox(0,0)[b]{$x^0_{2}$}}
\put(31,15){\makebox(0,0)[b]{$x^*_{1}$}}
\put(44,15){\makebox(0,0)[b]{$x^*_{2}$}}
\put(64,20){\makebox(0,0)[b]{$u_0$}}
\put(62,37){\makebox(0,0)[b]{$u^*$}}
\put (25,15){\vector(1,0){5}}
\put (25,15){\vector(0,1){5}}
\put (25,14){\makebox(0,0)[b]{$\bullet$}}
\put (52.7,15){\vector(-1,0){5}}
\put (52.7,15){\vector(0,1){5}}
\put (52.7,14){\makebox(0,0)[b]{$\bullet$}}

\put(38,40){\makebox(0,0)[b]{$u$}}
\put(85,10){\makebox(0,0)[b]{$x$}}

\put(7,10){\makebox(0,0)[b]{$-1$}}
\put(70,15){\line(0,1){2}}
\put(72,10){\makebox(0,0)[b]{$1$}}
\put(35,-5){\makebox(0,10)[b]{ Figure a). \; \rm Control of 
two points of sign change.}}
\end{picture}

\vspace{0.5cm}
\subsubsection*{Further results}
 In the following, we derive two results that generalize Theorem \ref{th1.1}.
\begin{cor} \label{th1.2} 
Let $ u_0,\,u^* \in  H_a^1 (-1,1).$ 
Assume that $ u_0$ and $u^*$ have finitely many points of sign change and
 the amount of points of sign change of $u^*$ is less than the one
  of $u_0.$
 Then, for any $ \eta > 0$ there exist $ T=T(\eta,u_0,u^*) > 0$ and  a piecewise static multiplicative control $ \alpha=\alpha(\eta,u_0,u^*) \in L^\infty (Q_T) $ such that  the solution $u$ to $\eqref{Psemilineare}$ satisfies
$$
\parallel u(\cdot, T) - u^* \parallel_{L^2 (-1,1)} \; \leq \; \eta.
$$
\end{cor}
{\bf Proof (of Corollary \ref{th1.2}).}
Corollary \ref{th1.2} easily follows from Theorem \ref{th1.1}. Indeed, all the target states described in Corollary \ref{th1.2} can be approximated in $L^2(-1,1)$ by those in Theorem \ref{th1.1}. $\;\;\diamond$\\

In the following Remark \ref{rem th1.2} we clarify the statement of Corollary \ref{th1.2}.
\begin{rem}\label{rem th1.2}\rm
 We note that by Corollary \ref{th1.2} we can steer the system \eqref{Psemilineare} from the initial state $u_0$ toward those states whose points of change of sign are organized in any order. 
We explain the statement of Corollary \ref{th1.2} by the following example. Let us denote by $x_l^0,\;l=1,\dots,n$ the points of sign change of $u_0.$ Let us consider an interval $ (-1, x^0_1)$ of positive values of $u_0$ followed by an interval $ (x^0_1, x^0_2)$ of negative values of $ u_0 (x)$, which in turn is followed by  an interval $ (x^0_2, x^0_3)$ of positive values of $ u_0 (x)$ and so forth. Then, the merging of the respective two points of sign change $ x^0_1 $ and  $ x^0_2$ will result in one single  interval $ (-1, x^0_3)$ of positive otherwise negative values. 
\end{rem}
In Figure $b)$ we describe one of the situations discussed in
 Remark \ref{rem th1.2} in the particular case $-1=x_0^0<x_1^0<x_2^0<x_3^0=1.$\\
 
\vspace{-1.9cm}
 
\setlength{\unitlength}{1mm}
\begin{picture}(150,55)(0,0)
\linethickness{1pt}


\put (2,15){\vector(1,0){85}}
\qbezier(10,19)(16,32)(23,18)
\qbezier(23,18)(34,1)(45,10)
\qbezier(45,10)(65,25)(69.5,37)
{
\qbezier(10,27)(35,4)(70,31)}
\put (10,05){\line(0,1){35}}
\put (40,05){\vector(0,1){35}}
\put (70,05){\line(0,1){35}}


\put(23,10){\makebox(0,0)[b]{$x^0_{1}$}}
\put(52,10){\makebox(0,0)[b]{$x^0_{2}$}}
\put(36,8){\makebox(0,0)[b]{$u_{0}$}}
\put(35,17){\makebox(0,0)[b]{$u^*$}}

\put(38,38){\makebox(0,0)[b]{$u$}}
\put(85,10){\makebox(0,0)[b]{$x$}}

\put(7,10){\makebox(0,0)[b]{$-1$}}
\put(70,15){\line(0,1){2}}
\put(72,10){\makebox(0,0)[b]{$1$}}
\put(42,0){\makebox(0,10)[b]{{ Figure b).} \;$u_0, {
u^*}$:\;{\rm 
merging of the points of change of sign.} }}
\end{picture}

 %

The following  approximate controllability property can be deduced from  Corollary \ref{th1.2}.
\begin{cor}\label{cor CK}
 Let $u_0$ and $u^*$ be given  in $ L^2 (-1,1)$.  Then, for any $ \eta > 0$ there exists $u^\eta_0 \in H_a^1 (-1,1)$ such that $\parallel u^\eta_0 - u_0 \parallel_{L^2 (-1,1)} \; <\eta,$ and there exist $ T=T(\eta,u_0,u^*)> 0$ and a piecewise static multiplicative control $\alpha=\alpha(\eta,u_0,u^*)\in L^\infty(Q_T)$ such that the solution $u$ to
$$
\left\{\begin{array}{l}
\displaystyle{u_t-(a(x) u_x)_x =\alpha(x,t)u+ f(x,t,u)\,\quad \mbox{ in } \; Q_T \,:=\,(-1,1)\times (0,T) }\\ [2.5ex]
\displaystyle{
B.C.
}\\ [2.5ex]
\displaystyle{u(0,x)=u^\eta_0 \in H_a^1 (-1,1)\,\qquad\qquad\qquad\qquad\quad\;\; \,x\in(-1,1)}~.
\end{array}\right.
  $$ 
satisfies
$$
\| u(\cdot, T) - u^* \|_{L^2 (-1,1)} \; \leq\eta.
$$
\end{cor}

{\bf Proof (of Corollary \ref{cor CK}).} The proof of Corollary \ref{cor CK} is similar to one of Corollary 2.2 of \cite{CFK}. 
$\qquad\qquad\qquad\qquad\qquad\qquad\qquad\qquad\qquad\qquad\qquad\qquad\qquad\qquad\qquad\qquad\qquad\qquad\diamond$


\subsection{Control strategy for the proof of the main result (Theorem \ref{th1.1})}\label{Sec Control Strategy} 
Let us consider the initial state $u_0\in H^1_{a}(-1,1)$ and the target state $u^*\in H^1_{a}(-1,1).$ Both data have  
$n$ points of sign change. 
Set $x^0_0:=-1$,  $x^0_{n+1}:=1,$ and consider the set of points of sign change of $u_0,$
 $X^0=(x^0_1,\ldots,x^0_n)$ where
$-1=x^0_0<x^0_l<x^0_{l+1}\leq  x^0_{n+1}=1$ for all $l=1,\ldots,n.$ 
Similarly, let $x^*_0:=-1$,  $x^*_{n+1}:=1,$ and
consider the set of target points 
$X^*=(x^*_1,\ldots,x^*_n),$ where
$-1=x^*_0<x^*_l<x^*_{l+1}\leq x^*_{n+1}=1$ for all $l=1,\ldots,n$. 
\subsubsection*{Some notations}
Let us introduce some notations.
\subsubsection*{Notation for space intervals}
\noindent Set 
$\displaystyle \rho^*_0=\min_{l=0,\ldots,n}\{x_{l+1}^*-x_{l}^*,\;x_{l+1}^0-x_{l}^0\},$
we define $a_{0}^*:=-1+\frac{\rho^*_0}{2}\; \text{ and } \; b_0^*:=1-\frac{\rho^*_0}{2},$
then $(a_{0}^*,b_{0}^*)\subset(-1,1).$

\subsubsection*{Notation for time intervals}
\noindent Given $N\in\N,$ 
for every $
(\tau_1,\ldots,\tau_N)=(\tau_k)_1^N\in\R^N_+,\;\;
(\s_1,\ldots,\s_N)=(\s_k)_1^N\in\R^N_+,  \,(
\footnote{ $\R^N_+=\{(a_1,\ldots,a_N)\,|\,a_k\in\R,\,a_k>0, \,k=1,\ldots,N\}$.})
$
we define 
\begin{equation}\label{N1}
T_0:=0,\qquad S_k:=T_{k-1}+\s_k, \qquad T_k:=S_k+\tau_k,\qquad\quad k=1,\ldots,N.  \qquad(\footnote{ We note that $T_0=0,\, S_1=\s_1,\, T_1=\s_1+\tau_1,$ and $\displaystyle S_k=\sum_{h=1}^{k-1}(\s_{h}+\tau_{h})+\s_k,\; T_k=\sum_{h=1}^{k}(\s_h+\tau_h),$ $\forall k=2,\ldots,N.$})
\end{equation} 
Noting that $0=T_0<S_{k}<T_{k}\leq T_N,\;k=1,\ldots,N,$ we consider 
the following partition of $[0,T_N]$ in $2 N$ intervals:
\begin{equation}\label{N2}
[0,T_N]=[0,S_1]\cup[S_1,T_1]\cup\cdots\cup[T_{N-1},S_N]\cup[S_{N},T_N]=\bigcup_{k=1}^N\left(\mO_k \cup \mE_k\right),
\end{equation}
where, for every $k=1,\ldots,N,$ 
we have set $\mO_k:=[T_{k-1},S_k]$ ($k^{th}$ odd interval)
and $\mE_k:=[S_k,T_k]$ ($k^{th}$ even interval). 
\subsubsection*{Notation for parabolic domains}
\noindent For every $k=1,\ldots,N,$ let us set $Q_{\mE_k}:=\displaystyle (-1,1)\times 
[S_k,T_k]$ and
 $Q^*_{\mE_k}:=\displaystyle (a_0^*,b_0^*)\times 
 [S_k,T_k]\subset 
 Q_{\mE_k}.$\\
 Let $Q_{\mE}:=\displaystyle (-1,1)\times {\small{\small \bigcup_{k=1}^N}}[S_k,T_k]$ and
 $Q^*_{\mE}:=\displaystyle (a_0^*,b_0^*)\times {\small{\small \bigcup_{k=1}^N}}[S_k,T_k]\subset 
 Q_{\mE}\,.$
\subsubsection*{Outline and main ideas for the proof of Theorem \ref{th1.1}}
The proof of Theorem \ref{th1.1} 
 uses the partition introduced in \eqref{N1}-\eqref{N2} and  two alternative control actions: 
on the even interval $\mE_k=[S_k,T_k]$ we choose suitable initial data, $w_k$, in pure diffusion problems ($v\equiv 0$) as control parameters to move the points of sign change to their desired location (see Section \ref{pure diffusion}), whereas on the odd interval $\mO_k=[T_{k-1},S_k]$ we give a smoothing result to preserve the reached points of sign change and attain such $w_k$'s as intermediate final conditions,
using piecewise static multiplicative controls $\alpha_k,\,\alpha_k\neq0$ (see Section \ref{sec smoth}).
The complete proof of Theorem \ref{th1.1} is achieved in Section \ref{proof main result}.
In the following figure we outline the iterative control strategy used to prove Theorem \ref{th1.1}, for simplicity, in the case of two points of sign change and Dirichlet boundary conditions, that is, in the $(WDP)$ case, with $\beta_1=\gamma_1=0.$\\ 

\vspace{0.5cm}
\begin{figure}[htbp]
\vspace{-1.3cm}
\begin{center}
\includegraphics[width=15.5cm, height=6.5cm
]{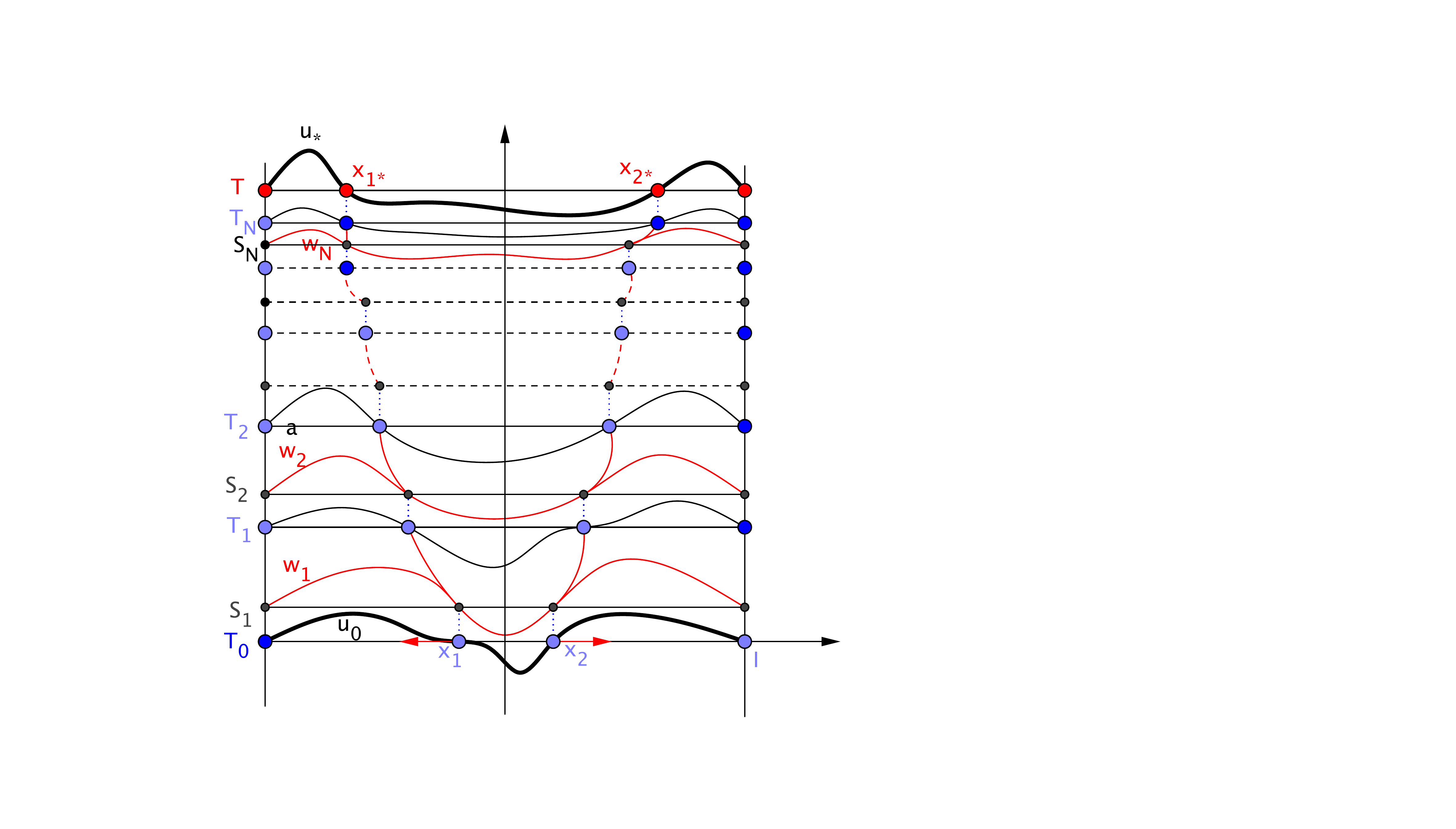}
\caption{ Iterative control strategy
}
\end{center}
\end{figure}

\vspace{-1.0cm}
 \subsection{Controllability for initial-value 
pure diffusion problems 
 }\label{pure diffusion} 
Let us fix a number $\beta\in (0,1)$ to
 be used in whole the paper.
Let $N\in\N.$ For any fixed $
(\s_1,\ldots,\s_N) 
\in \R^{N}_+,$ 
let us consider a generic $(\tau_1,\ldots,\tau_N)
\in \R^{N}_+$ and,
 recalling 
\eqref{N1}-\eqref{N2}, 
 let us introduce the following initial value 
pure diffusion problems on disjoint time intervals \\
\begin{equation}\label{2.1}
\left\{\begin{array}{l}
\displaystyle{
\quad   w_t \; = \; (a(x) w_x)_x
\;  + \; f(x,t,w), \;\;\;\;
 {\rm in} \;\;\;Q_{\mE
 } = 
 \displaystyle (-1,1)\times {\small{\small \bigcup_{k=1}^N}}[S_k,T_k]
}\\ [2.5ex]
\displaystyle{
\begin{cases}
\begin{cases}
\beta_0 w(-1,t)+\beta_1 a(-1)w_x(-1,t)= 0 \quad \;\:t\in 
\displaystyle {\small{\small \bigcup_{k=1}^N}}(S_k,T_k)\, \\
\qquad\qquad\qquad\qquad\qquad\qquad\qquad\qquad\qquad\qquad\qquad\qquad\quad\quad(\mbox{for }\, WDP)\\
\gamma_0\, w(1,t)\,+\,\gamma_1\, a(1)\,w_x(1,t)= 0
 \qquad\quad\, t\in \displaystyle {\small{\small \bigcup_{k=1}^N}}(S_k,T_k)
\end{cases}
\\
\quad a(x)w_x(x,t)|_{x=\pm 1} = 0\,\,\qquad\qquad\qquad\;\;\,\, t\in  \displaystyle {\small{\small \bigcup_{k=1}^N}}(S_k,T_k)
\;\;\quad\;\;\;\;\,(\mbox{for }\, SDP)
\end{cases}
}\\ [2.5ex]
\displaystyle{\quad w\mid_{t = S_k} \; = w_{k}(x),
\qquad\qquad\qquad\qquad\quad\;\; \:\,x\in(-1,1),\quad k = 1, \ldots, N}~.
\end{array}\right.
\end{equation}
Let us suppose that the assumptions $(A.3)$ and $(A.4)$ hold. Moreover,
we will consider the initial data $w_{k} $ and times 
 $\tau_k\; (\tau_k=T_k-S_k),
 $\, $ k = 1, \ldots, N,$ 
as control parameters, where the $ w_{k} $'s belong to $H^1_a(-1,1)\cap C^{2+ \beta} ([a_0^*,b_0^*]) \;(\footnote{ We recall the following spaces of H\"older continuous functions (see also \cite{FR}):
$$C^{\beta} ([a_0^*,b_0^*]):=\left\{w\in C ([a_0^*,b_0^*])\,:\:\sup_{x,y\in[a_0^*,b_0^*]}\frac{|w(x)-w(y)|}{|x-y|^\beta}
<+\infty\right\},$$
$$C^{2+\beta} ([a_0^*,b_0^*]):=\left\{w\in C^2 ([a_0^*,b_0^*])\,:\:w''\in C^{\beta} ([a_0^*,b_0^*])
\right\}.$$}),$ 
with 
$a_{0}^*=-1+\frac{\rho^*_0}{2}$ and $b_0^*=1-\frac{\rho^*_0}{2}$ $\;\displaystyle(\rho^*_0=\min_{l=0,\ldots,n}\{x_{l+1}^*-x_{l}^*,\;x_{l+1}^0-x_{l}^0\})$.
\begin{defn}\label{def ds}
We call 
solution of \eqref{2.1} the function defined in $\displaystyle (-1,1)\times {\small{\small \bigcup_{k=1}^N}}[S_k,T_k]
$ 
as
$$w(x,t)=W_k(x,t),\;\;\; \forall (x,t)\in(-1,1)
\times [S_{k}, T_{k}
],\;\;k=1,\ldots,N,$$
where $W_k$  is the unique strict solution on $(-1,1)\times [S_{k}, T_{k}
]$ of the $k^{th}$ problem in \eqref{2.1}.
\end{defn}
\begin{rem}\label{ds}\rm
We observe that a 
 solution of \eqref{2.1} is a collection of solutions of a finite number of
 problems which are set on disjoint time intervals. Therefore, it is 
 independent of the choice of $(\s_k)^N_1.$ We prefer to give the following Definition \ref{defTID} for a fixed $(\s_k)^N_1,$ just for technical purposes that will be clear in the sequel (see Theorem \ref{th2.2}). 
\end{rem}
\begin{defn}\label{defTID}
Let $u_0\in H^1_a(-1,1)$ be a function with the 
 $n$ points of sign change $x_l^0,\;l=1,\dots,n$. 
 For every fixed $N\in\N$ and $(\s_k)^N_1\in\R^N_+,$ we call a finite \lq\lq family of Times and Initial Data'' 
of 
\eqref{2.1} associated with $u_0,$ a set of the form
  $\left\{(\tau_k)_1^{N},(w_{k})_1^{N}\right\}
  $ 
  such that 
\begin{description}
\item[$\star$] \quad$(\tau_k)_1^N\in\R^N_+;$
\item[$\star$] \quad for all $k = 1, \ldots, N,$$\;\;w_{k}\in H^1_a(-1,1)\cap C^{2+ \beta} ([a_0^*,b_0^*]) 
$  satisfies the following:
\begin{enumerate}
\item 
$w_1$ and $u_0$ have the same points of sign change, in the same order as the points of sign change of $u_0;$
\item for $k = 2, \ldots, N,$ $ w_{k} (\cdot) $ and $w(\cdot, T_{k-1}
)$ have the same points as the points of sign change, in the same order of sign change of $u_0,$ where $w$ is the solution of \eqref{2.1} on $
\displaystyle (-1,1)\times {\small{\small \bigcup_{h=1}^{k-1}}}
[S_h,T_h]
.$
\end{enumerate}
\end{description}
\end{defn}
All Section \ref{ADCP} of this paper is devoted to the proof of the following Theorem \ref{th2.2}.
\begin{thm}\label{th2.2} 
Let $ 
{u}_{0}\in 
H^1_a(-1,1)$ have $n$ points of sign change at $ x^0_{l} \in (-1,1),\,  l=1,\ldots,n
$ 
with $$-1:=x^0_0<x^0_l<x^0_{l+1}\leq x^0_{n+1}:=1,\quad l=1,\ldots,n.$$
Let $ x^*_{l} \in (-1,1),\,  l=1,\ldots,n,$ be such that $-1:=x^*_0<x^*_l<x^*_{l+1}\leq x^*_{n+1}:=1.$ 
Then, for every $\ve>0$ there exist $N_{\ve}\in\N$ and a finite family of times and initial data 
$ \{(\tau_k)^{N_{\ve}}_1, (w_{k})^{N_{\ve}}_1\}$ 
 such that, for any $(\s_k)^{N_\ve}_1\in\R^N_+,$   
the 
solution $w^{\ve}$ of 
problem \eqref{2.1} 
 satisfies
$$ w^{\ve} (x, T_{N_{\ve}}
)=0 \qquad \Longleftrightarrow \qquad x=x_l^{\ve},\;\;l=1,\ldots,n,$$
for some points $x_l^{\ve}\in (-1,1)$,  with $-1:=x^{\ve}_0<x^{\ve}_l<x^{\ve}_{l+1}\leq x^{\ve}_{n+1}:=1$ for  $l=1,\ldots,n,$ such that
$$\displaystyle\sum_{l=1}^n|x^*_l-x^{\ve}_l|<\ve.$$
Moreover, $w^\ve(\cdot, T_{N_{\ve}}
)$ has the same order of sign change as $u_0.$ 
\end{thm}
In Section \ref{proof main result} we will use Theorem \ref{th2.2} to prove Theorem \ref{th1.1}.

\subsection{A control result 
to preserve the reached points of sign change and to obtain suitable smooth 
intermediate data 
}\label{sec smoth}
 In this section we introduce a smoothing result to preserve the reached points of sign change and attain smooth intermediate final conditions $w_k$'s.
  
Let $N\in\N.$ For any fixed $
(\tau_1,\ldots,\tau_N) 
\in \R^{N}_+,$ 
let us consider a generic $(\s_1,\ldots,\s_N)
\in \R^{N}_+$ and, for
$ k = 1, \ldots, N,$ recalling 
\eqref{N1}-\eqref{N2}, 
given $u_{k-1},\,r_{k-1}\in H^1_a(-1,1),\,\alpha_k\in  L^\infty((-1,1)\times  [T_{k-1},S_k]
),$ let us introduce the following problem
\begin{equation}\label{u pres}
\left\{\begin{array}{l}
\displaystyle{
\quad   u_t \; = \; 
(a(x) u_x)_x \; + \; \alpha_k(x,t)  u \; + \; f(x,t,u)
\;\;\;  {\rm in} \;\,Q_{\mO_k}= (-1, 1) \times [T_{k-1}, T_{k-1}+\s_k],  }\\ [2.5ex]
\displaystyle{
\begin{cases}
\begin{cases}
\beta_0 u(-1,t)+\beta_1 a(-1)u_x(-1,t)= 0 \quad \;\:t\in 
(T_{k-1}, T_{k-1}+\s_k)\, \\
\qquad\qquad\qquad\qquad\qquad\qquad\qquad\qquad\qquad\qquad\qquad\qquad\qquad\!\qquad(\mbox{for }\, WDP)\\
\gamma_0\, u(1,t)\,+\,\gamma_1\, a(1)\,u_x(1,t)= 0
 \qquad\quad\, t\in 
 (T_{k-1}, T_{k-1}+\s_k)
\end{cases}
\\
\quad a(x)u_x(x,t)|_{x=\pm 1} = 0\,\,\qquad\qquad\qquad\;\;\,\, t\in  
(T_{k-1}, T_{k-1}+\s_k)
\quad\quad(\mbox{for }\, SDP)
\end{cases}
}\\ [2.5ex]
\displaystyle{\quad u\mid_{t = T_{k-1}} \; = u_{k-1} + r_{k-1}\in H^1_a(-1,1)
}~,
\end{array}\right.
\end{equation}
where we recall that $S_k=T_{k-1}+\s_k,$ and we suppose that the assumptions $(A.3)$ and $(A.4)$ hold.
All Section \ref{nonneg gen} of this paper is devoted to the proof of the following Theorem \ref{th preserving}.
\begin{thm}\label{th preserving} 
Let $ u_{k-1}, r_{k-1}, w_k\in 
H^1_a(-1,1).$ Let $u_{k-1}$ and $w_{k}$ have the same $n$ points of sign change,
in the same order.
Then, for every $\eta>0$ there exist 
$\s_k>0$, $C_k\geq1$, and a piecewise static control $\alpha_k\in L^\infty((-1,1)\times(T_{k-1},S_k))$
(depending only on  $\eta, u_{k-1}$, and $w_k$)
such that
$$\|U_k (\cdot, S_k) - w_{k}(\cdot)\|_{L^2(-1,1)}\leq\eta+ C_k
\|r_{k-1}\|_{L^2(-1,1)}, 
$$ 
 where 
 $U_k$ is the solution of \eqref{u pres} on $(-1,1)\times[T_{k-1},S_k]$.
 \end{thm}
 In Section \ref{proof main result} we will use Theorem \ref{th preserving} to prove Theorem \ref{th1.1}.
 \begin{rem}\rm
 We note that Theorem \ref{th preserving}, in the particular case $r_{k-1}=0,$ gives the approximate controllability in the subspace of states with the same 
changes of sign, 
in the same order of sign change. 
\end{rem}
 \subsection{Proof of 
 Theorem \ref{th1.1}
 }\label{proof main result}

As soon as we prove Theorem \ref{th2.2} and Theorem \ref{th preserving}, combining these two results the proof of Theorem \ref{th1.1} is easily obtained, using an idea introduced in Section 3.3 of \cite{CFK}, through an intermediate result, Lemma \ref{pre proof th} (this lemma is similar to Lemma 3.1 of \cite{CFK}).
In this section we will avoid some repetitions, so we put only a sketch of the iterative idea and we invite the reader to see Section 3.3 of \cite{CFK}, that contains every technical detail. 
  \begin{lem}\label{pre proof th}
 Let 
$u_0\in H^1_a(-1,1)$ be a function with 
 $n$ points of sign change,
let $\displaystyle \left\{(\tau_k)_1^{N},(w_{k})_1^{N}\right\}$ be a finite family of times and initial data 
of 
\eqref{2.1} associated with $u_0,$ and
let\\ $\displaystyle w:(-1,1)\times\bigcup_{k=1}^N[S_k,T_k]
\longrightarrow\R$ be the 
solution of \eqref{2.1}.
Then for every $\delta>0$ there exists $\s_\delta=(\s_k)_1^N\in \R^N_+$ and $\alpha_\delta\in L^\infty((-1,1)\times(0,T_N))$
such that, denoting by $u_\delta:(0,1)\times[0,T_N]\rightarrow\R$ the solution of \eqref{Psemilineare} with bilinear control $\alpha_\delta,$ we have
 \begin{equation}\label{Lem1}
\|u_{\delta}(\cdot,T_{k})-w(\cdot, T_{k})\|_{L^2(-1,1)}
\leq\delta, \;\;\;	\forall k=1,\ldots,N.
\end{equation}
 \end{lem}
 {\bf Skech of the proof of Lemma \ref{pre proof th}.} 
 Fix $\displaystyle \left\{(\tau_k)_1^{N},(w_{k})_1^{N}\right\}$ and $\delta>0.$ 
 Let us consider the partition of $[0,T_N]$ in $2 N$ intervals introduced in 
 \eqref{N1}-\eqref{N2}.
 In particular, we will show that the bilinear control $\alpha_\delta$ has the following expression
 $$
 \alpha_\delta(x,t)=
  \begin{cases}
\quad   \alpha_k^\delta(x,t) \;\;\;\;\;\;\,
 {\rm in} \;\;\;Q_{\mO_k
 } = (-1,1) \times [T_{k-1},S_k] 
 , \;\:k=1,\ldots,N,
\\
\quad 0 
\qquad\quad\quad\;\;{\rm in} \,\;\;\;Q_{\mE_k
 } = (-1,1) \times [S_k,T_k],
  \;\:\;\;\;\;k=1,\ldots,N. 
\end{cases}
$$
In the following, for every $k=1,\ldots,N,$ we will consider the following problem on $Q_{\mE_k}$
$$
 \begin{cases}
\quad   u_t \; = \; (a(x)u_{x})_x \; 
 \; + \; f(x,t,u),
&\quad  {\rm in} \;\;\;Q_{\mE_k} = (-1,1) \times [S_k, T_k], \;\; 
\\
\;\;B.C.
\\
\quad u\:\mid_{t = S_k} \; = w_k+p_k, 
\end{cases}
$$
where $ w_k\in H^1_a(-1,1)\cap C^{2+ \beta} ([a_0^*,b_0^*]) ,\; p_k\in H^{1}_a(-1,1)$ are given functions,
and we will represent its 
solution 
as the sum of two functions $ w (x,t) $ and $h (x,t)$, which solve the following  problems in $Q_{\mE_k}$
\begin{equation}\label{lem h}
   \begin{cases}
\quad   w_t=(a(x)w_{x})_x
+f(x,t,w),
\\
\;\; \; B.C. 
\\
\quad w \:\mid_{t = S_k} \; = w_k,
\end{cases}
\!\!\!\!\!
   \begin{cases}
\quad   h_{t} \; = \; (a(x)h_{x})_x 
\; \!\!+\!\! \; (f(x,t,w+h) - f(x,t,w)),
\\
\;\; \; B.C.
\\
\quad h \:\mid_{t = S_k} \; = p_k. 
\end{cases}\!\!\!\!\!
\end{equation}
Multiplying by $ h$ the equation of the second problem of (\ref{lem h}) and integrating by parts over $ Q_{\mE_k},$ since $\left[a(x)h_x(x,t)h(x,t)\right]_{-1}^{1}\leq0,$ by \eqref{fsigni}
we obtain
\begin{align*}
\int_{-1}^1 h^2 (x,T_k) dx 
\leq&\int_{-1}^1 p_{k}^2 (x) dx
+2 \int_{S_k}^{T_k}  \int_{-1}^1 (f(x,t,w+h) -  f(x,t,w))h \,dx dt \nonumber\\\leq &
 \int_{-1}^1 p_{k}^2 (x) dx
+ 2\nu \int_{S_k}^{T_k}  \int_{-1}^1  h^2 dx dt,\quad t\in (S_k,T_k),
\end{align*}
thus applying Gr\"onwall's inequality
we deduce
\begin{equation}\label{hGron}
\parallel h(\cdot,T_k) \parallel_{L^2 (-1,1)} \; 
\leq  \; e^{\nu\widetilde{T}} \parallel  p_k \parallel_{L^2 (-1,1)},\quad \text{ with }\;\; \displaystyle\widetilde{T}:=\sum_{k=1}^{N}\tau_{k},\;\;\;\;\;  k=1,\ldots,N.
\end{equation}
Using the energy estimate \eqref{hGron} we can conclude this proof proceeding exactly with the same technical and iterative proof of Lemma 3.1 of \cite{CFK}. 
$\qquad\qquad\qquad\qquad\qquad\qquad\quad\qquad\qquad\qquad\qquad\qquad\qquad\quad\diamond$
\section{
Proof of Theorem \ref{th2.2} }\label{ADCP}  
\noindent In this section we refer to the notation introduced in Section \ref{pure diffusion}.
The plan of this section is as follows: 
\begin{description}
\item In Section \ref{pre lem}, we start by a regularity result for the problem  \eqref{2.1}, contained in Proposition \ref{solution of dtp}. By Lemma \ref{l1} we give the existence of suitable {\it initial data $w_k$'s}  to be used in the proof of Theorem~\ref{th2.2}.
By Lemma \ref{p1} we construct the $n$ {\it curves of sign change} associated with 
the $n$ initial points of sign change.
\item In Section \ref{OPS}, we construct a suitable particular family of {\it times and initial data},
that allows to move the $n$ initial points of sign change towards the $n$ target points of sign change. 
In this section, we also introduce
the definitions of {\it gap and target distance functional}. 
\item In Section \ref{Th3}, after obtaining Proposition \ref{P2}, we show how to steer the points of sign change of the solution
arbitrarily close to the target points. 
\end{description}
\subsection{Preliminary results}\label{pre lem}
Let us prove the following Proposition \ref{solution of dtp}.
\begin{prop}\label{solution of dtp}
Let $[a_0^*,b_0^*]\subset (-1,1),$ let $k=1,\ldots,N.$ 
If $ w_{k}\in H^1_a(-1,1)\cap~C^{2+ \beta} ([a_0^*,b_0^*]),$ 
then the $k^{th}$ initial-value problem in \eqref{2.1} has a unique strict solution $W_k(x,t)$ on $\overline{Q_{\mE_k}}$
 and
 $$W_k \in {{\cal{H}}(Q_{\mE_k})}\cap C^{2+\beta, 1+\beta/2} (\overline{Q^*_{\mE_k}}
 )\,.\;\;\; (\footnote{ Let $Q^*_{\mE_k}=
 (a_0^*,b_0^*)\times[S_k,T_k]$. We define 
 $${\cal{H}}(Q_{\mE_k}):=L^{2}(S_k,T_k;
H^2_a(-1,1))\cap H^{1}(S_k,T_k;L^2(-1,1))\cap C([S_k,T_k];H^{1}_a(-1,1)),
$$
  $$C^{\beta, \frac{\beta}{2}} (\overline{Q^*_{\mE_k}}):=\left\{u\in C (\overline{Q^*_{\mE_k}})\,:\:\sup_{x,y\in[a_0^*,b_0^*]}\frac{|u(x,t)-u(y,t)|}{|x-y|^\beta}+\sup_{t,s\in [S_k,T_k]}\frac{|u(x,t)-u(x,s)|}{|t-s|^\frac{\beta}{2}}<+\infty\right\},
 $$
 $$C^{2, 1} (\overline{Q^*_{\mE_k}}):=\left\{u:\overline{Q^*_{\mE_k}}\longrightarrow \R\,:\exists\:u_{xx},u_x,u_t\in C(\overline{Q^*_{\mE_k}
 }
 )
    \right\},
 $$  
  $$C^{2+\beta, 1+\frac{\beta}{2}} (\overline{Q^*_{\mE_k}}
  ):=\left\{u\in C^{2,1} (\overline{Q^*_{\mE_k}})\,:\:u_{xx},u_x,u_t\in C^{\beta, \frac{\beta}{2}} (\overline{Q^*_{\mE_k}})
    \right\}.
 $$ 
})
$$
     \end{prop}
{\bf Proof.
}
We note that our assumptions permit to apply a well known {\it interior regularity} result, contained in Section 5 and 6 of Chapter V in \cite{LSU} (in particular see Theorem 5.4, pp. 448-449, and Theorem 6.1, pp. 452-453). Indeed, on the domain $Q^*_{\mE_k}
=(a_0^*,b_0^*)\times[S_k,T_k]\subset\subset Q_{\mE_k}
$ the equation in \eqref{2.1} is uniformly parabolic, thus the unique strict solution of the $k^{th}$ problem in \eqref{2.1}, $W_k \in {{\cal{H}}(Q_{\mE_k})},$ is bounded on $Q^*_{\mE_k}$, therefore there exists a positive constant $M_k,$ such that 
$$|W_k(x,t)|\leq M_k,\qquad\mbox{ for a.e. } (x,t)\in Q^*_{\mE_k}.$$
Let us set
$$
f_{M_k}(x,t,w):=\begin{cases}
\qquad   f(x,t,w), \;\;\;\;\qquad
 {\rm if} \,\;\;|w|\leq M_k,
\\
\;\;
\quad  f(x,t,M_k), \;\;\;\qquad
 {\rm if} \;\;\;w> M_k,\\
\;
\quad  f(x,t,-M_k), \,\;\qquad
 {\rm if} \;\;\;w<- M_k.\\
%
\end{cases}
$$
Thus, by the inequalities 
 \eqref{fsigni} and Remark \ref{rem1} we deduce that, $\mbox{ for a.e. } (x,t)\in Q_T, \forall w_1,w_2\in \R$,
 $$
 \big|f_{M_k}(x,t,w_1)-f_{M_k}(x,t,w_2)\big|\leq
\nu
(1+|M_k|^{\vartheta-1}+|M_k|^{\vartheta-1})|u-v|=L(\vt, M_k)|u-v|,
$$
where $L(\vt, M_k):=\nu
(1+2|M_k|^{\vartheta-1}),$ and $\nu$ is the constant of $\eqref{fsigni}.$
Then, $w\longmapsto f_{M_k}(x,t,w)$ is a Lipschitz continuous function on $Q^*_{\mE_k}$, and we can apply the aforementioned
interior regularity result of \cite{LSU} to the problem
$$
 \begin{cases}
\quad   w_t \; = \; (a(x) w_x)_x
\;  + \; f_{M_k}(x,t,w), \;\;\;\;
 {\rm in} \;\;\;Q^*_{\mE_k
 } = (a_0^*,b_0^*) \times [S_{k}, T_{k}
 ],
\\
\quad w\mid_{t = S_k} \; = w_{k}(x),
\end{cases}
$$
thus, the unique solution $W_k$ belongs to $C^{2+\beta, 1+\beta/2} (\overline{Q^*_{\mE_k}}).$ $\qquad\qquad\qquad\qquad\qquad\qquad\quad\diamond$\\

By a simple exercise we can obtain the following.
\begin{lem}[Existence of suitable initial data $w_k$'s]\label{l1} 
Let $x_l\in (-1,1),\;l=1,\ldots,n,$ be such that 
$-1:=x_0<x_l
<x_{l+1}\leq x_{n+1}:=1,\;l=1,\ldots,n,.$
Let $
(\lambda_1,\ldots, \lambda_{n})\in~\R^{n}$, $
(\omega_1,\ldots, \omega_{n})\in~\R^{n}\;$ be such that $\lambda_l\in\{-1,1\},\,\omega_l\in\{-1, 0, 1\},\:l=1,\ldots,n,
$ and $\lambda_l\lambda_{l+1}~<~0,\,l=1\cdots,n-1.$
Let $\displaystyle \widetilde{\rho}=
\min_{l=0,\ldots,n}\left\{x_{l+1}-x_l
\right\},\,\tilde{a}:=-1+\frac{\tilde{\rho}}{2}$ and $\tilde{b}:=-1+\frac{\tilde{\rho}}{2}$. 
Then, there exists $w\in H^1_a(-1,1)\cap C^\infty([\tilde{a},\tilde{b}])$ such that
\begin{itemize}
\item[$\star$] $w(x)=0\;\;\Longleftrightarrow\;\; x=x_l,\; 
 l=1,\ldots, n;$
\item[$\star$] $w'(x_l)=\lambda_l,$\; $w''(x_l)=\omega_l,\;
l=1,\ldots, n;$
\item[$\star$] $\|w\|_{C^{m}([\tilde{a},\tilde{b}])}\leq C(m,\widetilde{\rho}), \;\forall m\in \N$.
\end{itemize}
\end{lem}
\begin{lem}[Construction of  the
curves of sign change]\label{p1}
Let $
(\lambda_1,\ldots, \lambda_{n})\in \R^{n}$ be 
such that $\lambda_l\in\{-1,1\},\, l=1,\ldots,n,\,$ and  $\lambda_{l}\,\lambda_{l+1}<0,\,l=1,\ldots,n-1. 
$
Let 
$\widetilde{\rho}>0$ and let $x_l\in (-1,1),\, l=1,\ldots,n,$ be such that 
$-1:=x_0
<x_l<x_{l+1}\leq x_{n+1}:=
1, \, l=1,\ldots,n,$ and 
$\displaystyle\min_{l=0,\ldots, n}\left\{x_{l+1}
-x_l
\right\}=\widetilde{\rho}.$
Let $a^*_0:=-1+\frac{\widetilde{\rho}}{2}$ and $b^*_0:=1-\frac{\widetilde{\rho}}{2}.$ Let $w_k\in H^1_a(-1,1)\cap C^{2+\beta}([a^*_0,b^*_0])$ be such that 
\begin{itemize}
\item[$\star$] $w_k(x)=0\;\;\Longleftrightarrow\;\; x=x_l,\; 
 l=1,\ldots, n;$
\item[$\star$] $w_k'(x_l)=\lambda_l,\; 
l=1,\ldots, n;$
\item[$\star$] $\|w_k\|_{C^{2+\beta}([a^*_0,b^*_0])}\leq c,$ for some positive constant $c=c(\widetilde{\rho}).$ 
\end{itemize}
Let $T>0$ and let $w$
be the solution of 
\begin{equation}\label{u}
\begin{cases}
\quad   w_t \; = \; \left(a(x)w_{x}\right)_x \;  + \; f(x,t,w) \;\;\;\;
 \qquad{\rm in} \;\;\quad Q_{T} = (-1,1) \times (0, T)
\\
\quad
B.C.\\ 
\quad w(x,0)
=w_k(x)\qquad\qquad\qquad\;\qquad\;\,\quad\;\;\;x\in(-1,1) \,.
\end{cases}
(\footnote{ 
For the existence, uniqueness and regularity of problem \eqref{u} see Proposition \ref{solution of dtp}. 
}) 
\end{equation}
Then, for every $\rho\in(0,\widetilde{\rho}]$ there exist $\widetilde{\tau}=\widetilde{\tau}({\rho}) >0$ and $M=M(\rho)>0$ such that, for each $l=1,\ldots,n,$ 
there exists a unique solution $\xi_l 
:[0,\widetilde{\tau}]\longrightarrow\R\, 
$  
of the initial-value problem
$$\;\;\;\displaystyle\begin{cases}
 \dot{\xi_l}
(t) \; 
=-\frac{a(\xi_l(t))w_{xx}(\xi_l(t),t)}{w_x(\xi_l(t),t)}-a'(\xi_l(t)),\;\;\;\; 
 t\in[0,\widetilde{\tau}],
\\
 \xi_l
(0)=x_l, 
\end{cases}
$$
that satisfies
\begin{itemize}
\item[$\bullet$] $w(\xi_l(t),t)=0,\;\qquad \forall t\in [0,\widetilde{\tau}];$
\item[$\bullet$] $\xi_l
\in C^{1+\frac{\beta}{2}}([0,\widetilde{\tau}])$\quad and \qquad $\|\xi_l\|_{C^{1+\frac{\beta}{2}}([0,\widetilde{\tau}])}\leq M;$
\item[$\bullet$] $\|\xi_l (\cdot) - x_l \|_{C([0,\widetilde{\tau}])} <\frac{\rho
}{2}.$
\end{itemize}
%
\end{lem}
\begin{defn}\label{scgs}
We call the functions $\xi_l:[0,\widetilde{\tau}]\longrightarrow \R,\; l=1,\ldots,n,$ given by Lemma \ref{p1}, Curves of Sign Change associated with the set of initial points of sign change $X
=(x_1,\ldots,x_n)$.
\end{defn}
\begin{rem}\label{separeted curves}\rm
As consequence of Lemma \ref{p1}, since $\|\xi_l (\cdot) - x_l \|_{C([0,\widetilde{\tau}])} <\frac{\rho}{2}\,$ for each $l=1,\ldots,n,$ 
 we have 
\begin{equation*}
a_0^*<\xi_{l}(t)<\xi_{l+1}(t)<
b_0^*,\;\;\;\forall t\in[0,\widetilde{\tau}], \;\forall l=1,\ldots,n-1.
\end{equation*}
Therefore, two adjacent curves of sign change don't intersect, so they remain separated.
\end{rem}
{\bf Proof (of Lemma \ref{p1}).}\; Let us fix $\rho\in(0,\widetilde{\rho}].$
Due to Proposition \ref{solution of dtp}, 
the solution $w$ of \eqref{u} is such that 
$$
w\in C^{2+\beta,1+\frac{\beta}{2}}(\overline{(a_{0}^*,b_0^*)\times(0,T)})\qquad \text{ and } \qquad  \| w \|_{C^{2+\beta,1+\frac{\beta}{2}}(
\overline{(a_{0}^*,b_0^*)\times(0,T)})}\leq
K,
$$ 
for some positive constant $K=K( \| w_k \|_{C^{2+\beta}([a_{0}^*,b_0^*])})$ 
(see (6.8)-(6.12) on pp. 451-452 in \cite{LSU}).
Thus, since $ \| w_k \|_{C^{2+\beta}([a_{0}^*,b_0^*])}\leq  c(\widetilde{\rho}),$ 
we have that
\begin{equation}\label{uest}
\| w \|_{C^{2+\beta,1+\frac{\beta}{2}}(\overline{Q}_T)} \; \leq
K( \| w_k \|_{C^{2+\beta}([a_{0}^*,b_0^*])})
\leq \; C, 
\end{equation}
for some positive constant $C=C(\widetilde{\rho}).$\\ 
{\it Existence and regularity of curves of sign change.}
 For any fixed $l=1,\ldots,n,$
 since $w_x(x_l,0)=\lambda_l\neq0$ and $w_x(x,t)$ is a continuous function in $(x_l,0)\in \overline{(a_{0}^*,b_0^*)\times(0,T)}
 ,$  there exist $\delta_l\in(0,\min\left\{\frac{1}{2C},\rho\right\})$(\footnote{ $C$ is the constant present in \eqref{uest}.}) and
  $T_l>0$ 
 such that $w_x(x,t
)\neq0,\;\forall (x,t)\in [x_l-\delta_l, x_l+\delta_l]\times[0,T_l].$\\
For every $l=1,\ldots,n,$
we consider the Cauchy problems 
\begin{equation}\label{ODE}
\begin{cases}
\quad   \dot{\xi_l}
(t) \; = \;-\frac{ w_{t}(\xi_l
(t),t)}{w_x(\xi_l
(t),t)},\;\;\;\;
 t>0,
\\
\quad \xi_l
(0)=x_l\;.
\end{cases}
\end{equation}
Let 
$\displaystyle \delta:=\min_{l=1,\ldots,n}\delta_l,$
we note that 
$G(x,t):=-\frac{w_{t}(x,t
)}{w_x(x,t
)}
$ is continuous on $[x_l-\delta, x_l+\delta]\times [0,T_l].$ 
Therefore, for every $l=1,\ldots,n,$ the problem \eqref{ODE} has a solution $\xi_l$ of class 
$C^1$ 
 on some interval $[0,\tau_l],$ with $0<\tau_l\leq T_l
.$ 
Let $\displaystyle \tau:=\min_{l=1,\ldots,n} \tau_l,$ since
$G
\in C^{\frac{\beta}{2}}([x_l-\delta, x_l+\delta]\times [0,T_l]),$ we conclude that $\xi_l\in C^{1+\frac{\beta}{2}}([0,\tau])$ and there exists $M=M(\rho)>0$ such that $\|\xi_l\|_{C^{1+\frac{\beta}{2}}([0,\widetilde{\tau}])}\leq M.$\\
Moreover, 
since by \eqref{ODE} we deduce that
$$\frac{d}{dt}w(\xi_l
(t),t)=w_{t}(\xi_l
(t),t)+w_{x}(\xi_l
(t),t)\dot{\xi_l}
(t)=0, \forall t\in [0,\tau],
\;
\text{ and }
\;
w(\xi
_l(0),0)=
w_k(x_l)=0,$$
we obtain that
$$w(\xi_l(t),t)=0,\;\;\forall t\in[0,\tau].$$ 
Furthermore, since by \eqref{Superlinearit}
 $f(w(\xi_l(t)),t)=0,\;\forall t\in[0,\tau],$
we also have that, for every $t\in [0,\tau],$
$$\displaystyle
\dot{\xi_l}
(t)=-\frac{ w_{t}(\xi_l
(t),t)}{w_x(\xi_l
(t),t)}
=-\frac{ \left(a(\xi_l
(t))w_{x}(\xi_l
(t),t)\right)_x}{w_x(\xi_l
(t),t)}=
-\left[ a'(\xi_l
(t))+\frac{a(\xi_l
(t))w_{xx}(\xi_l
(t),t)}{w_x(\xi_l
(t),t)}\right]
.$$
{\it Uniform estimates for the curves of sign change.}
For any fixed $l=1,\ldots,n,$ we consider
the
uniform time 
$\widetilde{\tau}
 := \min\Big\{\Big(\frac{1}{2C}-\delta\Big)^\frac{2}{\beta}, \frac{\delta^2}{3}, \tau
 \Big\},
$ 
 $\widetilde{\tau}=\widetilde{\tau}(\rho)$ (independent of $l$) since $\delta<\min\left\{\frac{1}{2C},\rho\right\}.$ 
Recall that the function $t\mapsto w_x(x,t)$ belongs to $C^{\frac{\beta}{2}}([0,\widetilde{\tau}])$ and the function $x\mapsto w_x(x,t)$ belongs to $C^{1+\beta}([0,\widetilde{\tau}]).$ Thus, for every $(x,t)\in (x_l-\delta,x_l+\delta)
\times (0,\widetilde{\tau}),$ by \eqref{uest} we have
\begin{align}\label{uxest}
| w_{x}(x,t) - \lambda_l | &
\leq | w_{x}(x,t) - w_{x}(x,0)|  + | w_{x}(x,0) - w_{x}(x_l,0)|\nonumber\\
&\leq \|w\|_{C^{2+\beta,1+\frac{\beta}{2}}(\overline{Q}_T)}( t^{\frac{\beta}{2}}+|x - x_l|)\leq 
%
 C ( \widetilde{\tau}^{\frac{\beta}{2}}+\delta ).
\end{align}
Since
 $\delta<\frac{1}{2C}$ and
 $\widetilde{\tau}\leq \big(\frac{1}{2C}-\delta\big)^\frac{2}{\beta},$
we deduce
$C (\widetilde{\tau}^{\frac{\beta}{2}} + \delta) \leq \frac{1}{2},
$
so by \eqref{uxest} we obtain
$$\Big||w_{x}(x,t)| - |\lambda_l|\Big|\leq| w_{x}(x,t) - \lambda_l |
\leq\frac{1}{2}.$$
Therefore, for every $l = 1, \ldots, n,$ having in mind that $ | \lambda_l | = 1,$ we have
\begin{equation}\label{0.5}
|w_{x}(x,t)|\geq | \lambda_l |-\frac{1}{2}=\frac{1}{2}, \;\;\;\;\;\;\forall (x,t)\in (x_l-\delta,x_l+\delta)\times (0,\widetilde{\tau})\,.
\end{equation}
Then, by \eqref{uest} and \eqref{0.5}, keeping in mind that $ \widetilde{\tau}\leq\delta^2/3$ and $\delta<\min\left\{\frac{1}{2C},\rho\right\},$ for every $t\in [0,\widetilde{\tau}],$ we deduce 
\begin{equation}\label{delta1}
\hspace{-0.1cm}\;
|\xi_l (t) - x_l |=\left| \int_0^t\dot{\xi_l} (s)\,ds \right|   \leq  \int_0^{\widetilde{\tau}} \frac{ | w_{xx}(\xi_l(s),s) | }{|w_x(\xi_l(s),s)| } ds \; 
\leq  \;2
\widetilde{\tau}  \|w\|_{C^{2+\beta,1+\frac{\beta}{2}}(\overline{Q}_T)}
\!\! \leq\,
 2\widetilde{\tau} C<\frac{\widetilde{\tau}}{\delta}
 < \frac{\rho}{3}.
\end{equation}  
{\it Uniqueness \lq\lq a posteriori'' of the curves of sign change.}
We observe that, although one cannot claim uniqueness for the Cauchy problem \eqref{ODE}, 
{\it a posteriori} the $\xi_l$'s 
 turn out to be uniquely determined.
 Indeed, let $\bar{a}\in(-1,a_0^*]\;\text {and} \;\bar{b}\in[b_0^*,1),$
 setting $\xi_0(t)\equiv \bar{a},\;\xi_{n+1}(t)\equiv \bar{b},\,\forall t\in~[0,\widetilde{\tau}],$ 
 since by 
 \eqref{delta1} and Remark \ref{separeted curves} 
 one can apply the strong maximum principle for uniformly parabolic equations 
 on the domains 
 $\left\{(x,t)|x\in\left[\xi_l(t), \xi_{l+1}(t)\right],\,t\in[0,\widetilde{\tau})\right\},$ 
 for every $l=0,\ldots,n.$ The fact that the initial datum $w_k(x)$ doesn't change sign on $(x_l,x_{l+1})$ 
   implies that (thanks to the fact that $\bar{a}$ and $\bar{b}$ are arbitrary), for every $t^*\in [0,\widetilde{\tau}),$
$$w(x,t^*)=0\quad\Longleftrightarrow \quad x=\xi_l(t^*),\; l=0,\ldots,n+1,$$
 completing the proof of  Lemma \ref{p1}. 
$\qquad\qquad\qquad\qquad\qquad\qquad\qquad\qquad\qquad\qquad\qquad\diamond$
\subsection{Construction of Order Processing Steering sets}\label{OPS}
In the following we define {\it Order Processing Steering Times and Initial Data}
that permit to move the points of sign change towards the desired targets. 
In this section we use the notation introduced in Section \ref{Sec Control Strategy} 
and in Section \ref{pure diffusion}.\\ 

Given the initial state $u_0\in H^1_a(-1,1)$, 
let us consider the $n$ points of sign change of $u_0,$
 $X^0=(x^0_1,\ldots,x^0_n),$ where
$-1:=x^0_0<x^0_l<x^0_{l+1}\leq  x^0_{n+1}:=1,\,$ for $l=1,\ldots,n.$ Let us define, for every $l=1,\ldots,n,$
\begin{equation}\label{lambda}
\lambda(x^0_l)=
\begin{cases}
 1,\;\;\;\;
 \text{ if }  u_0(x)>0\; \text{ on } (x^0_l,x^0_{l+1}),
\\
-1,\;\, \text{ if } u_0(x)<0\; \text{ on } (x^0_l,x^0_{l+1}).
\end{cases}
\end{equation} 
Since $x_l^0,\,l=1,\ldots,n,$ are points of sign change, we note that $\lambda(x^0_{l+1})=-\lambda(x^0_{l}),\,l=1,\ldots,n-1.$\\

Let us set $x^*_0:=-1$ and $x^*_{n+1}:=1,$ and let us consider the set of $n$ target points 
$X^*=(x^*_1,\ldots,x^*_n),$ where
$-1=x^*_0<x^*_l<x^*_{l+1}\leq x^*_{n+1}=1,\;\; l=1,\ldots,n.$\\ 
Let $\displaystyle \rho^*_0=
\min_{l=0,\ldots,n}\{x^0_{l+1}-x^0_l,\, x^*_{l+1}-x^*_l
\},$ then let us set $\displaystyle a_{0}^*:=-1+\frac{\rho^*_0}{2}\;$ and $\;b_0^*:=1-\frac{\rho^*_0}{2}.$
\subsubsection*{Order Processing Steering Times and Initial Data} 

Let $\beta\in (0,1)$ be the number that was fixed at the beginning of Section \ref{pure diffusion}.  
Let $N\in\N$ and let us fix $(\s_k)^N_1\in\R^N_+.$ 
Now, we will construct 
a {\it finite
a family of Times and Initial Data}
for \eqref{2.1} associated with $u_0,$
  $\left\{(\tau_k)_1^{N},(w_{k})_1^{N}\right\}$ (see Definition \ref{defTID}), 
in order to move the points of sign change towards the desired targets. We denote by  ${\cal{W^*}}(u_0)$ the subclass of
 such special families, which we call  {\it Order Processing Steering Times and Initial Data} 
associated with $u_0$ and $X^*$.\\
In the following, we will define the times $S_k,T_k,\; k=1,\ldots,N,$ in the same way of 
\eqref{N1}.

 \begin{description}
 \item {\it Construction of $\big\{\tau_1, w_1\big\}$}.
 
By Lemma \ref{l1}, there exists $w_1\in H^1_a(-1,1)\cap C^{2+\beta}([a^*_{0},b_0^*])
,$ with $\|w_1\|_{C^{2+\beta}([a^*_{0},b_0^*])}\leq c_1,$ for some positive constant $c_1= c(
{\rho_0^*}
),$
 such that
\begin{itemize}
\item[$\star$] $w_1(x)=0\quad\Longleftrightarrow\quad x=x^0_l, \quad 
l=1,\ldots, n;$
\item[$\star$] $w_1'(x^0_l)=\lambda(x^0_l)a(x^0_l),\quad 
w_1''(x^0_l)=-\lambda(x^0_l)
\left[\mu_1
(x_l^*-x_l^0)
+a'(x^0_l)\right],
\;\;l=1,\ldots, n,
$\\
where
$
\mu_1
(x_l^*-x_l^0)=sgn(x_l^*-x_l^0)=
\begin{cases}
 1,\;\;\;\;
 \text{ if } \,x_l^0<x_l^*,\;\;\;\;\\
0,\;\;\;\;\text{ if }\,x_l^0=x_l^*,\\
-1,\,\;\text{ if }\,x_l^0>x_l^*.
\end{cases}
$
\end{itemize}
Let $w
$ be the solution to
$$
\begin{cases}
\quad   w_t \; = \; \left(a(x)w_{x}\right)_x \;  + \; f(x,t,w) \;\;\;\;
 \qquad
 (x,t)\in(-1,1) \times (S_1, +\infty)
\\
\quad \displaystyle{
B.C. } 
 \\
\quad w(x,S_1)
=w_k(x),\qquad\qquad\;\,\quad\;\;
\end{cases}
$$
where $S_1=\s_1.$
By Lemma \ref{p1}, for every $\displaystyle \rho\in(0,\rho_0^*)
$ there exist $\displaystyle\widetilde{\tau}_1=\widetilde{\tau}_1({\rho})>0,$ $M_1=M_1(
{\rho
})>0$ 
and
$n$  curves of sign change (associated 
to the points of sign change $X^{0}=(x^{0}_1,\ldots,x^{0}_n)$)
  $\xi^1_l
\in C^{1+\frac{\beta}{2}}([S_1,\widetilde{T}_1]), 
\, l=1,\ldots,n,$ with $\widetilde{T}_1=S_1+\widetilde{\tau}_1,$
such that $\|\xi^1_l
\|_{C^{1+\frac{\beta}{2}}([S_1,\widetilde{T}_1])}\leq M_{1},$ $w(\xi^1_l
(t),t)= 0\;\,\forall t\in [S_1,\widetilde{T}_1],$ and
\begin{equation}\label{xi}
\!\!\!\!\!
\begin{cases}
  \dot{\xi^1_l}
(t) \; = \;
-\left[ a'(\xi^1_l
(t))+\frac{a(\xi^1_l
(t))w_{xx}(\xi^1_l
(t),t)}{w_x(\xi^1_l
(t),t)}\right],\,
 t\in[S_1,\widetilde{T}_1],
 \vspace{.1cm}
\\
%
\xi_l^1
(S_1)=x^0_l\,.
\end{cases}
\end{equation}
Let us set
$\xi^1_0(t)\equiv 
a^*_0\;$ and $\;\xi^1_{n+1}(t)\equiv 
b^*_0$ on $[S_1,\widetilde{T}_1],$ by Remark \ref{separeted curves}, for every $l=1,\ldots,n-1,$ we have
\begin{equation}\label{sccs1}
a^*_0=\xi^1
_0(t)<\xi^1
_l(t)<\xi^1
_{l+1}(t)<\xi^1
_{n+1}(t)=b^*_0,\; \forall t\in [S_1,\widetilde{T}_1].
\end{equation}
Let us introduce the {\it Inactive Set} (\footnote{ The {\it  Inactive Set} $L^{0}_{IS}$ is the set of the indexes such that the corresponding points of sign change don't need to be moved.})
$$L^{0}_{IS}:=\{l\,|\,l\in\{1,\ldots,n\},\,x_l^0=x_l^*\}$$
and let us consider the set of the {\it stopping times}  
$$\Theta_{1}
:=\{s\in(0,\widetilde{\tau}_{1}]\,|\,\xi_{l}^1(S_1+s)=x_l^*,\;\text{ for some }\,l\in\{1,\ldots,n\}\backslash L^{0}_{IS}\}.
$$
Let us set
\begin{equation}\label{tau***1}
\displaystyle
\tau_{1}=\begin{cases}
\; \widetilde{\tau}_1
\qquad\qquad\quad \text{ if }\quad 
\;\;\Theta_{1}=\varnothing,
\\
\;\displaystyle\min \Theta_{1},
\qquad\qquad\text{ otherwise },
\end{cases}
\end{equation}
by \eqref{N1} we have
$T_{1}=S_1+\tau_1.$
 \item {\it Construction of $\big\{\tau_k, w_k\big\},\;k=2,\ldots,N$.}
 
By the previous step 
 we have obtained the vector $X^{k-1}=(x_1^{k-1},\ldots,x_n^{k-1}),$  where
$x^{k-1}_l:=\xi_l^{k-1}(T_{k-1})$ for $l=1,\ldots,n,$ and  
 $\xi_l^{k-1},$ defined on $[S_{k-1},T_{k-1}],$ are the $n$ {\it curves of sign change} associated with the initial state $w_{k-1}$ and to the set of points of sign change $X^{k-2}=(x^{k-2}_1,\ldots,x^{k-2}_n)$. Let us set $x^{k-1}_0:=a^*_0$,  $x^{k-1}_{n+1}=:b_0^*$, and $ \rho^*_{k-1}
=
\displaystyle\min_{l=0,\ldots,n}\{x^{k-1}_{l+1}-x^{k-1}_l,\,x^{*}_{l+1}-x^{*}_l
\}$.
Let us introduce the 
{\it Inactive Set} 
  $$L^{k-1}_{IS}
:=\{l,\;l\in\{1,\ldots,n\}|\exists h_l\in\{1,\ldots,k-1\}:
x^{h_l}_l=x^*_l\},\; 
$$
which consists of the indexes of the points of sign change that have already reached the corresponding target points $X^*=(x_1^*,\ldots,x_n^*)$ at some previous time (\footnote{ We note that $L^{k-1}_{IS}\subseteq\{1,\ldots,n\}$ 
 is an increasing family of sets.}),
so these points don't need to be moved.
Then, let us set
$$
\displaystyle
\mu_k
(x_l^*-x_l^0)=
\begin{cases}
 0\;\;\;\;
\;\;\quad\qquad\qquad\quad\quad\;\; \text{ if } \,l\in L_{IS}^{k-1}
,
 \;\;\;\;\\
 sgn(x_l^*-x_l^0)
\;\;\qquad\;\;\quad\text{ if }\, l\not\in L_{IS}^{k-1}\,.
\end{cases}
(\footnote{
We remark that the definition of $\mu_k, \,k=2,\ldots,N,$ is consistent with the one of $\mu_1.$})
$$
By Lemma \ref{l1} we can choose $w_k\in H^1_a(-1,1)\cap C^{2+\beta}([a_0^*,b_0^*]),$ with $\|w_k\|_{C^{2+\beta}([a_0^*,b_0^*])}\leq c_k,$
 for some positive constant $c_k= c(
{\rho^*_{k-1}}
),$ 
such that
\begin{itemize}
\item[$\star$] $w_k(x)=0\quad\Longleftrightarrow\quad x=x^{k-1}_l,\quad
  l=1,\ldots, n;$ 
\item[$\star$] $w_k'(x^{k-1}_l)=\lambda(x^{0}_l)a(x_l^{k-1}),$\;
 $w_{k}''(x^{k-1}_l)=\!\!-\lambda(x^0_l)
 \!\!\left[\mu_k(x_l^*
 -x_l^0)+a'(x^{k-1}_l)\right],
\,l=1,\ldots, n.$
%
\end{itemize}

Let $w
$ be the solution to
$$
\begin{cases}
\quad   w_t \; = \; \left(a(x)w_{x}\right)_x \;  + \; f(x,t,w) \;\;\;\;
 \qquad
 (x,t)\in(-1,1) \times (S_k, +\infty)
\\
\quad \displaystyle{B.C.
}
 \\
\quad w(x,S_k)
=w_k(x).
\end{cases}
$$
where $\displaystyle S_k=\sum_{h=1}^{k-1}(\s_h+\tau_h)+\s_k.$
By Lemma \ref{p1}, for every $\rho\in(0,\rho^*_{k-1}]$ there exist $
\widetilde{\tau}_k=\widetilde{\tau}_k(\rho)>0,$ $M_k=M_k(\rho)>0$ 
and $n$ curves of sign change 
(
associated 
to the points of sign change
$X^{k-1}=(x_1^{k-1},\ldots,x_n^{k-1})$), $\xi^k_l
\in C^{1+\frac{\beta}{2}}([S_k,T_{k}]),\, l=1,\ldots,n,$
with $\displaystyle 
\widetilde{T}_k=S_k+\widetilde{\tau}_k,$ 
such that $\|\xi^k_l
\|_{C^{1+\frac{\beta}{2}}([S_k,\widetilde{T}_{k}])}\leq M_{k},$ $w
(\xi^k_l
(t),t)=0 \;\forall t\in [S_k,\widetilde{T}_k],$ and 
\begin{equation}\label{xik}
\!\!\!\!\!\begin{cases}
\quad   \dot{\xi^k_l}
(t) \; = -\left[ a'(\xi^k_l
(t))+\frac{a(\xi^k_l
(t))w_{xx}(\xi^k_l
(t),t)}{w_x(\xi^k_l
(t),t)}\right],\;\;\;\;\\
\quad \xi_l^k
(S_k)=x^{k-1}_l\,.
\end{cases}
\end{equation}
Let us set
$\xi^k_0(t)\equiv 
a_0^*\;$  and  $\;\xi^k_{n+1}(t)\equiv 
b_0^*$ on $[S_k,\widetilde{T}_{k}],$ by Remark \ref{separeted curves}, 
for every $l=1,\ldots,n-1,$ 
we have that
 \begin{equation}\label{sccsk}
 a_0^*=\xi^k
_0(t)<\xi^k
_l(t)<\xi^k
_{l+1}(t)<\xi^k
_{n+1}(t)=b_0^*,\; \qquad\forall t\in [S_k,\widetilde{T}_k].
\end{equation}
Let us consider the set of the {\it stopping times} 
$$\Theta_{k}:=\{s\in(0,\widetilde{\tau}_{k}]\,|\,\xi_{l}^k(S_k+s)=x_l^*,\;\text{ for some }\,l\in\{1,\ldots,n\}\backslash L^{k-1}_{IS}\},$$
and let us set 
\begin{equation}\label{tau***k}
\tau_k=\begin{cases}
\;\widetilde{\tau}_{k}
\qquad\qquad\;\; \text{ if }\quad 
\;\;\Theta_{k}=\varnothing,
\\
\;\displaystyle\min \Theta_{k},
\qquad\qquad\text{ otherwise },
\end{cases}
\end{equation}
by \eqref{N1} we have $T_k=S_k+\tau_k\,.$
\end{description}
\begin{rem}\label{remaining}\rm
We note that $\tau_k<\widetilde{\tau}_{k}$ for at most $n$ values of $k\in\{1,\ldots,N\}.$
\end{rem}
{\it An important remark.} Let us give an important remark about ${\cal{W^*}}(u_0),$ that is, the previous subclass of
 special families, which we have called  {\it Order Processing Steering Times and Initial Data} 
associated with $u_0$ and $X^*,$ and we will show that a generic $\left\{(\tau_k)_1^{N},(w_{k})_1^{N}\right\}\in{\cal{W^*}}(u_0)$ moves the points of sign change towards the desired targets. 
\begin{rem}\label{rem steer}
\rm
We note that, for each index $l\not\in L_{IS}^{k-1}\,(k=1,\ldots,N),$ 
by \eqref{xi} and \eqref{xik} and the choice of the initial data $w_k$ we deduce that
  $$ \dot{\xi^k_l}
(S_k) \; 
=sgn(x^{*}_l-x^0_l),\;\;\;\;\xi_l^k
(S_k)=x^{k-1}_l.$$
If $x^0_l<x^*_l$, we have that $\dot{\xi^k_l}
(S_k) \; =1
>0,$ thus the initial conditions $w_k$ permit to move the points of sign change $x^{k-1}_l$ to the right towards $x^*_l.$  Similarly, if $x^*_l<x^0_l$, the initial condition $w_k$ permits to move the points of sign change to the left. 
\end{rem}
\subsubsection*{Curves of Sign Change, Gap and Target Distance functional}

Given $W^N= \{(\tau_k)_1^{N},(w_{k})_1^{N}\}
\in{\cal{W^*}}(u_0),$ we introduce the $n$ curves of sign change 
associated with $W^N$
as the functions $\;\displaystyle\xi_l^{W}:\bigcup_{k=1}^N
[S_k,T_k]\longrightarrow\R, \;l=1,\ldots,n,$ such that 
$$\displaystyle \xi_l^{W}(t)=
\xi^{k
}_l
(t
),\qquad S_{k}\leq t\leq T_k,\;\; k=1,\ldots,N,$$
where the {curves}
$\xi^{k
}_l
$ are previously been constructed. We also  set 
$\xi^{W}_0(t)\equiv 
0$ and $\xi^{W}_{n+1}(t)\equiv 
1.$
Moreover, by \eqref{sccs1} and \eqref{sccsk}, for $l=1,\ldots,n-1,$
we deduce that
$$a_0^*=\xi_0^{W}(t)<\xi_l^{W}(t)<\xi_{l+1}^{W}(t)< \xi_{n+1}^{W}(t)=b_0^*,\;\;\; \text{ for all } \displaystyle t\in \bigcup_{k=1}^N [S_k,T_k]
.$$
\begin{defn}
For all $W^N=
 \left\{(\tau_k)_1^{N},(w_{k})_1^{N}\right\}
\in{\cal{W^*}}(u_0)$ we define the \textbf{gap functional} by
$$\rho(W^N)=\min_{l=0,\ldots,n}\min_{t\in
\displaystyle  \cup_{k=1}^N[S_k,T_k]}\{\xi_{l+1}^{W}(t)-\xi_l^{W}(t)\}$$
and the \textbf{target distance 
functional} by
$$\displaystyle J^*(W^N)=
\sum_{l=1}^n |\xi_l^{W}(T_N)-x_l^*|\,.$$
\end{defn}
\subsection{Proof of Theorem \ref{th2.2} completed}\label{Th3}
Let us consider the initial state $u_0\in H^1_a(-1,1),$ 
and consider the set of points of sign change of $u_0,$
 $X^0=(x^0_1,\ldots,x^0_n)$ where
$0=x^0_0<x^0_l<x^0_{l+1}\leq  x^0_{n+1}=1$ for all $l=1,\ldots,n.$ 
Similarly, 
let consider the set of target points 
$X^*=(x^*_1,\ldots,x^*_n),$ where
$-1=x^*_0<x^*_l<x^*_{l+1}\leq x^*_{n+1}=1$ for all $l=1,\ldots,n$. 
Set 
$$\rho^*_0=\min_{l=0,\ldots,n}\{x_{l+1}^*-x_{l}^*,\;x_{l+1}^0-x_{l}^0\}$$
and let $\tau^*_0=\tau(\frac{\rho^*_0}{2})>0$ and $M^*_0=M(\frac{\rho^*_0}{2}
)>0$ 
be the positive time and constant of Lemma \ref{p1}, associated with $\rho=\frac{\rho^*_0}{2}
$. 
The following proposition is crucial to obtain the  proof of Theorem~\ref{th2.2}.
\begin{prop}\label{P2}
There exists $\ve^*_0\in (0,1)$ such that for all $\ve\in (0,\ve_0^*]$ 
and  $N\in \N,\,N>n$ 
there exists $W^N= \{(\tau_k)_1^{N},(w_{k})_1^{N}\}
\in{\cal{W^*}}(u_0)$  such that
$\rho(W^N)\geq \frac{\rho^*_0}{2}
$ and
 \begin{equation}\label{functional J}
J^*(W^N)
\leq \sum_{l=1}^{n}\,|x^0_{l}-x^*_{l}|+c_1(\ve)\,\sum_{k=1}^{N}\frac{1}{{k}^{1+\frac{\beta}{2}}}
-c_2(\ve)\sum_{k=n+1}^{N}\frac{1}{k}, 
\end{equation}
where $c_1(\ve)=
\frac{\ve\rho_0^*n}{4
\,s_\beta}$, $c_2(\ve)=
(\frac{\ve\rho_0^*}{4M^*_0
\,s_\beta})^{\frac{2}{2+\beta}}$, and  $\displaystyle 
 s_\beta=\sum_{k=1}^\infty\,\frac{1}{k^{
 1+\frac{\beta}{2}}}$.\\
Moreover, for  such $W^N,$ for $k=1,\ldots,N,$ we have that
\begin{equation}\label{W^N}
  \tau_k 
 \leq  "\widetilde{\tau}_{k}:= \left(\frac{\ve \rho_0^*}{4M^*_0\,s_\beta}\right)^{\frac{2}{2+\beta}}\frac{1}{k},
 \qquad (\footnote{ $\widetilde{\tau}_k$ and $\tau_k$ are defined in \eqref{tau***1} and \eqref{tau***k}, see also Remark \ref{remaining}.})
 \end{equation}
 and
 if
$l\in
L^{k-1}_{IS}$ the following inequality holds
\begin{equation}\label{stab}
|\xi^h_{l}(t)-x^*_{l}|\leq \ve\frac{\rho_0^*}{4},\quad \quad\forall t\in[S_{h},T_h], \quad \forall h=k,\ldots,N.\,
\end{equation}
\end{prop}
\begin{rem}\label{remInactive}\rm
We note that 
for each inactive index $l\in L_{IS}^{k-1} \,(k=1,\ldots,N)$ 
we have chosen the initial data
such that 
$\dot{\xi^h_l}
(S_h)\displaystyle
=0,\;\;
h=k,\ldots,N.$ So, by the inequality \eqref{stab} of Lemma \ref{P2}, the corresponding points of sign change remain forever near the target points that they have already reached. 
\end{rem}

We omit the proof of Proposition~\ref{P2}, because using Remark \ref{rem steer} and Remark \ref{remInactive} we can repeat a proof similar to that of Proposition 4.1 of \cite{CFK}.

We give the following definition.
\begin{defn}
A set of times and initial data 
 $W^N\in {\cal{W^*}}(u_0)$ is said to be {\em separating} if $\rho(W^N)\geq \frac{\rho^*_0}{2}$.
We set
${\cal{W_{S}^*}}(u_0):=\{W^N \in {\cal{W^*}}(u_0):\;\rho(W^N)\geq \frac{\rho^*_0}{2}\}.$
\end{defn}

Finally, we can prove Theorem~\ref{th2.2}.\\
{\bf Proof (of Theorem~\ref{th2.2}).}
We will prove  that  
\begin{equation*}
\forall \ve>0\;\;\exists\,N_{\varepsilon}\in\N\;\;\exists W^{N_\varepsilon}\in {\cal{W_{S}^*}}(u_0)
\;\;\text{such that}\;\; J^*(W^{N_\varepsilon})<\varepsilon
\;\;\; \text{ and } \;\;\; L^{N_{\ve}}_{IS}=\{1,\ldots,n\},
\end{equation*}
which implies the conclusion of Theorem~\ref{th2.2}. Arguing
by contradiction, suppose
$$\exists\, \ve>0\;\exists\,j\in \{1,\ldots,n\}:
\forall N\in\N\,,\;\forall W^{N}\in {\cal{W_{S}^*}}(u_0) \;\;\text {we have}\;\;
 |\xi^{N}_{j}(T_N)-x^*_{j}|>\ve.$$
Moreover, we can assume $\ve\leq\ve^*_0,$ 
where $\ve_0^*\in (0,1)$ is given by Proposition \ref{P2}. For eve\-ry $N>n,$ by Proposition \ref{P2}, there exists $W^N= \{(\tau_k)_1^{N},(w_{k})_1^{N}\}
\in{\cal{W_S^*}}(u_0)
$ such that
$\displaystyle\tau_k
 \leq~\widetilde{\tau}_{k}=~\Big(\frac{\ve \rho_0^*}{4M^*_0\,s_\beta}\Big)^{\frac{2}{2+\beta}}\frac{1}{k},$
 $\,k=1,\ldots,N$ and, by \eqref{functional J},
we obtain 
$$ \ve< |\xi^N_{j}(T_{N})-x^*_{j}|\leq J^*(W^N)\leq \sum_{l=1}^{n}\,|x^0_{l}-x^*_{l}|+c_1\,\sum_{k=1}^{N}\frac{1}{{k}^{1+\frac{\beta}{2}}}-c_2\,\sum_{k=n+1}^{N}\frac{1}{k} .
$$
Since $\displaystyle\,\sum_{k=1}^{\infty}\frac{1}{k}=+\infty\,,$  the previous inequality gives a contradiction. 
\hfill $\diamond$




\section{ Proof of Theorem \ref{th preserving} 
}\label{nonneg gen}
This section is devoted to the proof of Theorem \ref{th preserving}, obtained in Section \ref{regular section} after proving Lemma \ref{lem preserving}. Let us start with Section \ref{AppA}, where we recall some preliminary results obtained in \cite{F1} (for $(SDP)$) and in  \cite{F2} (for $(WDP)$). In this section we use the notation introduced in Section \ref{Wp}.

\subsection{
Some spectral properties and some estimates for the 
semilinear degenerate problem \eqref{Psemilineare}}\label{AppA}
Let us observe that the semilinear problem $(\ref{Psemilineare})$  can be recast 
as
$$
 \left\{\begin{array}{l}
\displaystyle{u^\prime(t)=A\,u(t)+\phi(u)\,,\qquad  t>0 }\\ [2.5ex]
\displaystyle{u(0)=u_0\in H^1_a(-1,1),
}~
\end{array}\right.
$$
where the operator $(A,D(A))$ is defined in \eqref{DA} and, for every $u\in {{\cal{H}}(Q_T)},
$ 
$$
\phi(u)(x,t):=f(x,t,u(x,t))\;\;\;\forall (x,t)\in Q_T.
$$
Let us consider the operator $(A_0,D(A_0))$ defined as 
\begin{equation}\label{DA0}
\left\{\begin{array}{l}
\displaystyle{
D(A_0)=D(A)
}\\ [2.5ex]
\displaystyle{A_0\,u=(au_x)_x, \,
\,\, \forall \,u \in D(A_0)}~,
\end{array}\right.
\end{equation}
for this operator 
the following Proposition \ref{str cont} and Proposition \ref{spectrum} is obtained in \cite{CFproceedings1} for $(SDP)$ (\footnote{ In the $(SDP)$ case, in \cite{CFproceedings1}, 
we showed that for this result it is sufficient that the diffusion coefficient
  $a(\cdot)$ satisfies 
  the assumption 
   $(A.4_{SD})$ with $\xi_a(x)
   =\int_0^x\frac{ds}{a(s)}
\in L^{1}(-1,1)$ instead of 
 $\xi_a
    \in L^{2\vartheta-1
    }(-1,1).$
}) and in \cite{CF2} for $(WDP).$

\begin{prop}\label{str cont}
$(A_0, D(A_0))$ is a closed, self-adjoint, dissipative operator with dense domain in $L^2 (-1,1)$.
Therefore, $A_0$ is the infinitesimal generator of a strongly continuous semigroup 
of bounded linear operator on $L^2(-1,1)$.
\end{prop}
Proposition \ref{str cont} permits to obtain the following.
\begin{prop}\label{spectrum}
There exists an increasing sequence $\{\lambda_p\}_{p\in\N},$ with
$\lambda_p\longrightarrow +\infty \mbox{ as } p \, \rightarrow\infty\,,$
such that the eigenvalues of the operator $(A_0,D(A_0))$
are given by $\{-\lambda_p\}_{p\in\N}$, and the corresponding eigenfunctions $\{\omega_p\}_{p\in\N}$ form a complete orthonormal system in $L^2(-1,1)$.
\end{prop}
\begin{rem}\label{Legendre}\rm
In the case 
$a(x)=1-x^2$, that is in the case of the Budyko-Sellers model, 
 the orthonormal eigenfunctions of the operator $(A_0,D(A_0))$ are reduced to Legendre's polynomials $Q_p(x)$, and the eigenvalues are $\mu_p=(p-1)p, \, p\in\N.$ $Q_p(x)$ is equal to $\sqrt{\frac{2}{2p-1}} L_p(x),$ where $L_p(x)$ is assigned by \textit{Rodrigues's
formula}:
$L_p(x)=\frac{1}{2^{p-1} (p-1)!} \frac{d}{dx^{p-1}} \, (x^2-1)^{p-1}, \; p \geq 1.$
\end{rem}

\n 
{Some estimates.} By the next Lemma \ref{sob3} and Lemma \ref{esist glob} (obtained in \cite{F1} and \cite{F2})
we can deduce the Proposition \ref{f in L2}. 
\begin{lem}\label{sob3}
Let $T>0$ and $\vartheta\geq1.$ Let $a \in C([-1,1])\cap C^1(-1,1)$ such that the assumption $(A.4)$ ($(A.4_{SD})$ or $(A.4_{WD})$) holds, 
then
$
{\cal{H}}(Q_T)\subset L^{2\vartheta}(Q_T)
$
\, and 
$$
\|u\|_{L^{2\vartheta}(Q_T)}
\leq c\,
 T^{\frac{1}{2\vartheta}}\,\|u\|_{{\cal{H}}(Q_T)},
$$
where $c$ is a positive constant.
\end{lem}
\begin{lem}\label{esist glob}
Let $T>0,\; u_0\in H^1_a(-1,1)$ and let $\alpha\in L^\infty(-1,1).$ 
The strict solution $u\in {\cal{H}}(Q_T)$ of system \eqref{Psemilineare}, under the assumptions $(A.3)-(A.4)$,
satisfies the following estimate
$$\|u\|_{{\cal{H}}(Q_{T})}\leq 
c\,e^{k T}\|u_0\|_{1,a},$$
where
$c=c(\|u_0\|_{1,a})
$ and \;
 $k$ are positive constants.
\end{lem}
So, we can deduce the following.
\begin{prop}\label{f in L2}
Let $T>0, $
\,$u_0\in H^1_a(-1,1)$ and let $\alpha\in L^\infty(-1,1).$
Let
$u\in {\cal{H}}(Q_T)$ the strict solution $u\in {\cal{H}}(Q_T)$ of system \eqref{Psemilineare}, under the assumptions $(A.3)-(A.4)$. 
Then, the function
$(t,x)\longmapsto f(t,x,u(t,x))$
belongs to $L^2(Q_T)$ and the following estimate holds
$$\int_{Q_T}|f(x,t,u(x,t))|^2\,dx\,dt\leq Ce^{2k \vt T}
T\|u_{0}\|^{2\vt}_{1,a
}\;, 
%
$$
where $C=C(\|u_{0}\|_{1,a}) \text{ and } k$ are positive constants.
\end{prop}
{\bf Proof.} Using Lemma \ref{sob3} and Lemma \ref{esist glob} we obtain
$$
\int_0^T  \int_{-1}^1   f^2(x,t,w) dx dt\leq \gamma_*^2 \int_0^T  \int_{-1}^1   |w|^{2\vt} dx dt\\
\leq 
c
T\|w\|^{2\vt}_{{{\cal{H}}(Q_T)}}
\leq 
 Ce^{2k \vt T}
T\|u_{0}\|^{2\vt}_{1,a
}, 
$$
where $c=c(\|u_{0}\|_{1,a}),\,C=C(\|u_{0}\|_{1,a}) \text{ and } k$ are positive constants. $\qquad\qquad\qquad\quad\quad\quad\diamond$

\subsection{
Proof 
}\label{regular section}
In this section 
we reformulate the problem \eqref{u pres}, using a lighter notation than one introduced in Section \ref{sec smoth} in the statement of Theorem \ref{th preserving}, in the following way  
\begin{equation}\label{u pres proof}
   \begin{cases}
\quad   u_t \; = \; \left(a(x)u_{x}\right)_x \; + \; \alpha(x,t)  u \; + \; f(x,t,u)
&\quad  {\rm in} \;\;\;Q_{T}= (-1, 1) \times (0, T),  
\\
\quad B.C. 
&\quad t \in (0, T),
\\
\quad u\mid_{t = 0} \; = u_{in} + r_{in},
\end{cases}
\end{equation}
where $(0, T)$ is a generic time interval, 
 $u_{in},\,r_{in}\in 
H^1_a(-1,1),$ and $u_{in}$ has $n$ points of sign change at $ x_{l} \in (-1,1),\,l=1,\ldots,n, 
$ 
with $-1:=x_0<x_l<x_{l+1}\leq x_{n+1}:=1.$ Moreover, we will denote the target state by $ \overline{u} \in H^1_a(-1,1)$ instead of 
$ w_k.$ 

Throughout this section, we represent the solution $u(x,t)$ 
of \eqref{u pres proof} as the sum of two functions $ w (x,t) $ and $h (x,t)$, which solve the following problems in $Q_T:$
\begin{equation}\label{w}
 \!\!\! \begin{cases}
\;   w_t=\left(a(x)w_{x}\right)_x +  
\alpha
 w+f(x,t,w)
\\
\; B.C. 
\\
\; w\mid_{t = 0} \; = u_{in}, 
\end{cases}
\!\!\!\!\!   \begin{cases}
\;   h_{t}=\left(a(x)h_{x}\right)_x+\alpha h+(f(x,t,w+h) - f(x,t,w))
\\
\; B.C. 
\\
\; h\mid_{t = 0} \; = r_{in}. 
\end{cases}\!\!\!\!\!\!\!\!
\end{equation}

Let us start with the following Lemma \ref{lem preserving}.
\begin{lem}\label{lem preserving} 
Let $\overline{u}\in 
H^1_a(-1,1)$ have the same points of sign change as $u_{in}$ in the same order of sign change. Let us suppose that
\begin{equation}\label{ass}
\exists \delta^*>0:\;\delta^*\leq\frac{\overline{u}(x)}{u_{in}(x)}< 1, \qquad  
\forall \,x\in \left(-1,1\right)\backslash \bigcup_{l=1}^n\left\{x_l\right\}. 
\end{equation}
Then, for every $\eta>0$ there exist a small time $T=T(\eta,u_{in}, \overline{u})>0$ 
and a static bilinear control $\tilde{\alpha}(x),\,\tilde{\alpha}=\tilde{\alpha}(\eta,u_{in}, \overline{u})
\in C^2([-1,1])
$ such that
\begin{equation}\label{ineq pres}
\|u (\cdot, T) - \overline{u}(\cdot)\|_{L^2(-1,1)}\leq \eta+\sqrt{2}\|r_{in}\|_{L^2(-1,1)},
\end{equation}
 where $u$ is the corresponding solution of \eqref{u pres proof} on $Q_T.$
\end{lem}

\n
{\bf Proof.}
Let us represent the solution $ u $ of \eqref{u pres proof} 
as the sum of two functions $ w (x,t) $ and $h (x,t)$, which solve the two problems introduced in \eqref{w}, respectively.
For this proof it need to obtain some preliminary estimates that will deduced in the Step. 1.\\
{\bf Step 1:} {\it Evaluation of  $\| w \|_{C ([0,T]; L^2 (-1,1))}$ and  
 $\| h (\cdot, T) \|_{L^2 (-1,1)}$.}
   For every $\alpha\in L^\infty (Q_T),$ with $\alpha(x,t)\leq0\;\; \forall (x,t)\in Q_T,$
multiplying by $ w $ the equation in the first problem of (\ref{w}) 
and integrating by parts, using \eqref{fsigni}, 
since $\left[a(x)w_x(x,t)w(x,t)\right]_{-1}^{1}\leq0$
 we obtain 
\begin{align*}
\frac{1}{2}\int_0^t \int_{-1}^1 (w^2)_t dx ds &=
\int_0^t \int_{-1}^1 \left(a(x)w_{x}\right)_xw\, dx ds+
 \int_0^t \int_{-1}^1 \alpha w^2 dx ds+\int_0^t  \int_{-1}^1   f(x,s,w)w dx ds
\\
&\leq 
-\int_0^t \int_{-1}^1 a(x)w^2_{x}\, dx ds+
\nu\int_0^T  \int_{-1}^1   w^2 dx dt\leq \nu
\int_0^T  \int_{-1}^1   w^2 dx dt,
\end{align*}
where $ \nu $  is the constant of (\ref{fsigni}). Then, for $T\in (0, \frac{1}{4\nu})$ we deduce
\begin{multline*}\label{4.8}
\!\!\!\!\!\!\int_{-1}^1 w^2 (x,t) dx \leq \int_{-1}^1 u_{in}^2 (x) dx + 
2\nu\int_0^T  \int_{-1}^1   w^2 dx dt\\
\!\leq \int_{-1}^1 u_{in}^2 (x) dx  + 2\nu
T \|w\|^2_{C([0,T],L^2(-1,1))}\leq  \parallel  u_{in} \parallel^2_{L^2 (-1,1)}+ \frac{1}{2} \|w\|^2_{C([0,T],L^2(-1,1))},\;\;t\in(0,T),
\end{multline*}
so,
\begin{equation}\label{4.10}
\parallel w \parallel_{C ([0,T]; L^2 (-1,1))} \; \leq \; \sqrt{2} \parallel  u_{in} \parallel_{L^2 (-1,1)}.
\end{equation}

Proceeding as for the estimate \eqref{4.10} we evaluate $ \| h (\cdot, T) \|_{L^2 (-1,1)}.$ Namely, multiplying by $ h$ the equation of the second problem of (\ref{w}) and integrating by parts over $ Q_T,$ using \eqref{fsigni},
for $T\in (0,\frac{1}{4\nu})$ and for every $t\in(0,T),$ since $\left[a(x)h_x(x,t)h(x,t)\right]_{-1}^{1}\leq0$ we obtain
\begin{align*}
\int_{-1}^1 h^2 (x,t) dx 
\leq&\int_{-1}^1 r_{in}^2 (x) dx
+2 \int_0^T  \int_{-1}^1 (f(x,t,w+h) -  f(x,t,w))h \,dx dt \nonumber\\\leq &
 \int_{-1}^1 r_{in}^2 (x) dx 
+ 2\nu \int_0^T  \int_{-1}^1  h^2 dx dt 
\leq 
 \|r_{in} \|^2_{L^2(-1,1) } 
+2\nu T \|h\|^2_{C([0,T];L^2(-1,1))}
\nonumber\\
\leq &
\|r_{in} \|^2_{L^2(-1,1) } 
+\frac{1}{2}  \|h\|^2_{C([0,T];L^2(-1,1))}.
\end{align*}
 Hence,  
\begin{equation}\label{h}
 \| h (\cdot, T) \|_{L^2 (-1,1)}\leq \parallel h \parallel_{C ([0,T]; L^2 (-1,1))} \; \leq \; \sqrt{2} \parallel  r_{in} \parallel_{L^2 (-1,1)}.
\end{equation}
\noindent
{\bf Step 2:} {\it Choice of the bilinear control $\alpha$.}\;
Let us consider
the following function defined on $[-1,1]$
\begin{equation*}
\alpha_0 (x) \; = \;  \left\{ \begin{array}{ll}
\log \left( \frac{\overline{u} (x)}{u_{in} (x) } \right), \;\;\;\; &  {\rm for } \;\; x \neq  -1, 1, x_l  \; \;(l = 1\ldots, n)   \\ 
0, \;\;\;\; &  {\rm for} \;\; x =  -1, 1, x_l\; \;(l = 1\ldots, n) \,.
\end{array}
\right.  
\end{equation*} 
Using the assumption \eqref{ass}, we deduce that 
$\alpha_0\in L^\infty (-1,1)$  
and
$
\alpha_0 (x)\leq 0,
\text{ for every } x\in [-1,1].$ 
Now, we select the following bilinear control
$$
\alpha (x,t) \; := \; \frac{1}{T} \alpha_0 (x)\,.
$$

For every fixed $x\in (-1,1),$ by the classical technique for solving first order ODEs, applied to the equation $w_t(x,t)=\frac{\alpha_0(x)}{T}w(x,t)+\left(\left(a(x)w_{x}(x,t)\right)_x +  
f(x,t,w)\right)\;t\in (0,T),$
we compute at time $T$ the solution $w$ to the first problem in \eqref{w}, so the following representation formula holds for every $x\in(-1,1)$
\begin{equation}\label{WT}
w (x,T) \; = \; 
e^{\alpha_0(x)}u_{in}(x) \; + \int_0^T e^{\alpha_0 (x) \frac{(T - \tau)}{T}} \big((a(x)w_{x})_x (x, \tau) + f (x,\tau, w(x,\tau))\big) d\tau
.
\end{equation}
Let us show that 
$w(\cdot, T) \rightarrow \overline{u} $ in $L^2 (-1,1),$ as $ T \rightarrow 0^+.$
In advance, since $ \alpha_0 (x)\leq 0$ let us note that by the above formula we deduce 
\begin{multline}\label{uT}
\|w (\cdot,T)-\overline{u} (\cdot)\|^2_{L^2(-1,1)}=\int_{-1}^1 \left(\int_0^T e^{\alpha_0 (x)
\frac{(T - \tau)}{T}} ((a(x)w_{x})_x (x, \tau) + f (x,\tau,w(x,\tau))) d \tau \right)^2 dx\\ \; \leq \;
T \parallel  \left(a(\cdot)w_{x}\right)_x  + f (\cdot,\cdot,w) \parallel^2_{L^2 (Q_T)}.
\end{multline}
Let us prove that the right-hand side of (\ref{uT}) tends to zero 
 as $ T \rightarrow 0^+$.\\
{\bf Step 3:} {\it Evaluation of $\|\left(a(\cdot)w_{x}\right)_x\|^2_{L^2 (Q_T)}\!.$} 
Let us suppose, without loss of generality, 
that $\alpha_0$ satisfies the further properties:
\begin{equation}\label{Ass alpha}
 \alpha_0\in C^2 ([-1,1]) \text{  with }\displaystyle \lim_{x\rightarrow\pm1}\frac{\alpha_0'(x)}{a(x)}=0\;\; \text{ and }\;\;\lim_{x\rightarrow\pm1}\alpha_0'(x)a'(x)=0 \;\;
 . 
\end{equation} 
Moreover, let us consider the $(WDP)$ problem with $\beta_0\gamma_0\neq0$ in the assumption $({A.4}_{WD}).$ These assumptions will be removed in Step.4.

Multiplying by $\left(a(x)w_{x}\right)_x$ the equation in the first problem in \eqref{w}, with $ \alpha(x,t) = \frac{\alpha_0(x)}{T} \leq 0,$ integrating 
over $ Q_T,$ 
we have
\begin{multline}\label{a_x}
\parallel  \left(a(\cdot)w_{x}\right)_x  \parallel^2_{L^2 (Q_T)}
\leq \int_0^T \int_{-1}^1 w_t \left(a(x)w_{x}\right)_x dx dt \; - \;
\frac{1}{T} \int_0^T \int_{-1}^1 \alpha_0 w \left(a(x)w_{x}\right)_x dx dt \\+ \: 
\frac{1}{2}  \int_0^T  \int_{-1}^1   f^2(x,t,w) dx dt
+ \frac{1}{2}  \int_0^T  \int_{-1}^1  \left|\left(a(x)w_{x}\right)_x\right|^2 dx dt\,.
\end{multline}
Let us estimate the first two terms of the right-hand side of \eqref{a_x}. Integrating by parts and using the sign condition $\beta_0\beta_1\leq0$ and $\gamma_0\gamma_1\geq0$ 
we deduce
\begin{align}\label{deg2}
\!\!\!\!\!\int_0^T \int_{-1}^1 w_t \left(a(x)w_{x}\right)_x &dx dt =\int_0^T \left[w_t \left(a(x)w_{x}\right)\right]_{-1}^1 dt- \frac{1}{2}\int_0^T \int_{-1}^1 a(x)\left(w^2_{x}\right)_t dx dt \nonumber\\
&\leq\frac{1}{2}\frac{\gamma_1}{\gamma_0}a^2(1)u^2_{in x}(1)-\frac{1}{2}\frac{\beta_1}{\beta_0}a^2(-1)u^2_{in x}(-1)+ \frac{1}{2} \int_{-1}^1 a(x)u^2_{in x} dx dt,
\end{align}
moreover, 
by \eqref{Ass alpha} and \eqref{4.10} we obtain
\begin{align}\label{deg3}
 \int_0^T \int_{-1}^1 &\alpha_0(x) w \left(a(x)w_{x}\right)_x dx dt = 
 -\int_0^T \int_{-1}^1 \alpha_0(x)a(x)w^2_{x} dx dt 
  -\frac{1}{2}\int_0^T \int_{-1}^1 \alpha'_0(x)a(x)\left(w^2\right)_{x} dx dt \nonumber\\
&\geq -\frac{1}{2}\int_0^T \int_{-1}^1 \left(\alpha'_0(x)a(x)w^2\right)_x dx dt+ \frac{1}{2}\int_0^T \int_{-1}^1 \left(\alpha''_0(x)a(x)+\alpha'_0(x)a'(x)\right)w^2 dx dt \nonumber\\
&\geq-\frac{1}{2}\sup_{x\in [-1,1]}\left|\alpha''_0(x)a(x)+\alpha'_0(x)a'(x)\right|  \int_0^T \int_{-1}^1w^2 dx dt\nonumber\\
 &\geq-T\sup_{x\in [-1,1]}\left|\alpha''_0(x)a(x)+\alpha'_0(x)a'(x)\right|
\|u_{in}\|^2_{L^2 (-1,1)}\,.
\end{align}
From \eqref{uT}, applying \eqref{a_x}-\eqref{deg3} and Proposition \ref{f in L2}, we deduce
\begin{align}\label{uT2}
\|w (x,T)&-\overline{u} (x)\|^2_{L^2(-1,1)}\leq
2T\left( \| \left(a(\cdot)w_{x}\right)_x\|^2_{L^2 (Q_T)}  + \|f (\cdot,\cdot,w) \|^2_{L^2 (Q_T)}\right)\nonumber\\
&\leq 2T\left(\frac{\gamma_1}{\gamma_0}a^2(1)u^2_{in x}(1)-\frac{\beta_1}{\beta_0}a^2(-1)u^2_{in x}(-1)+
|u_{in}|^2_{1,a}\right)\nonumber\\
&+4T\left(\sup_{x\in [-1,1]}\left|\alpha''_0(x)a(x)+\alpha'_0(x)a'(x)\right| \|u_{in}\|^2_{L^2 (-1,1)}+ Ce^{2k \vt T}
\|u_{in}\|^{2\vt}_{1,a
}\right),
\end{align}
where $C=C(\|u_{in}\|_{1,a}) \text{ and } k$ are the positive constants of Proposition \ref{f in L2}.\\
\noindent {\bf Step 4:} {\it Convergence of $ w (\cdot, T) $ to $\overline{u}(\cdot)$.} 
The previous estimates of Step.3 hold also in the $(WDP)$ case with $\beta_0\gamma_0=0$ and in the $(SDP)$ case. In effect, the simple weighted Neumann boundary condition permits some simplifications in the previous estimates (see \eqref{deg2}), and in particular in the last estimate \eqref{uT2}.\\
In order to remove the assumption \eqref{Ass alpha}
we observe that
we can approximate in $ L^2 (-1,1)$ the reaction coefficient $\alpha_0\in L^\infty(-1,1),$ introduced in Step.2, by a sequence of 
uniformly bounded functions $\displaystyle \{ \alpha_{0j}\}_{j \in \N}\subset C^2( [-1, 1])$ such that 
$$\alpha_{0j} (x) \leq 0\;\; \forall x\in [-1,1],\; \;\alpha_{0j} (\pm1)=0,\,\displaystyle \lim_{x\rightarrow\pm1}\frac{\alpha_{0j}^\prime(x)}{a(x)}=0\;\; \text{ and }\;\;\lim_{x\rightarrow\pm1}\alpha_{0j}^\prime(x)a'(x)=0.$$ 
We remark that, for every $j\in \N$, the representation formula \eqref{WT} still holds for the corresponding solutions $w_j:$
$$w_j (x,T) \; = \; 
e^{\alpha_{0j}(x)}u_{in}(x) \; + \int_0^T e^{\alpha_{0j} (x) \frac{(T - \tau)}{T}} \big((a(x)w_j{_{x}})_x (x, \tau) + f (x,\tau, w_j(x,\tau))\big) d\tau.
$$
Let us fix $\eta>0,$ then, making use of the following limit relation
$$
\bar{u}_j(x):=e^{\alpha_{0j} (x)} u_{in} (x) 
 \; \longrightarrow  \;
e^{\alpha_0 (x)} u_{in}(x)
\; = \; \overline{u} (x) \;\;{\rm in} \;\; L^2 (-1,1) \;\; {\rm as} \; j \rightarrow \infty,
$$
there exists $j^*\in\N$ and a positive constant $K=K(\alpha_{0j^*},\alpha'_{0j^*},\alpha''_{0j^*}, u_{in}),$ given by \eqref{uT2} (written with $\alpha_{0j^*}$ instead of $\alpha_{0}$), such that we deduce that
\begin{multline}\label{w(T)}
\!\!\!\!\!\!\|w_{j^*} (x,T)-\overline{u} (x)\|_{L^2(-1,1)}\leq \|w_{j^*}(x,T)-\overline{u}_{j^*} (x)\|_{L^2(-1,1)}+\|\overline{u}_{j^*} (x)-\overline{u} (x)\|_{L^2(-1,1)}\\
\leq KT+\frac{\eta}{2}
\leq \eta,  
\end{multline}
for $T\leq \frac{\eta}{2K}.$   
Keeping in mind that $u_{j^*}=w_{j^*}+h_{j^*},$ by combining 
\eqref{w(T)} and \eqref{h} we obtain that there exist a small time $T=T(\eta,u_{in}, \overline{u})<\min\left\{ \frac{\eta}{2K}, \frac{1}{4\nu}\right\}$ 
and a static bilinear control $\tilde{\alpha}=\tilde{\alpha}(\eta,u_{in}, \overline{u})
\in C^2([-1,1])
\left(\tilde{\alpha}(x,t)=\frac{\alpha_{j^*}(x)}{T},\;\forall (x,t)\in Q_T\right)$ such that we obtain the conclusion. 
\hfill$\diamond$\\

Now, using the notation in \eqref{u pres proof}, let us give the proof of
Theorem \ref{th preserving}.\\

\noindent {\bf Proof of Theorem \ref{th preserving}.} Let us fix $\eta>0$. For any $\displaystyle \rho\in\left[0,\frac{\rho_0}{2}\right),\,$ with $\displaystyle \rho_0:= 
\!\!\min_{l=0,\ldots,n} \big\{ x_{l+1}-x_{l}\big\},$ let us consider
the  set 
$\displaystyle \Omega_{\rho}:=
\bigcup_{l=0}^{n} (x_l + \rho, x_{l+1} -\rho).$ 
Since $\overline{u},u_{in}\in H^1_a(-1,1),$ there exist
\\$\overline{u}_\eta, u^\eta_{in}\in C^1([-1,1])$ 
 such that 
 \begin{itemize}
 \item[$\star$]  $
 {\overline{u}_\eta}(x)=0 \iff x=x_l,\;\;l=0,\ldots,n+1\;\;\text{ and }\;\;\parallel \overline{u}_{\eta}-\overline{u} \parallel_{L^2 (-1,1)}<\frac{\eta}{4};$
\item[$\star$] 
${u_{in}^\eta}(x)=0 \iff  x=x_l,\; l=0,\ldots,n+1,
\; \;|({u_{in}^\eta})^\prime(x_l)|=1,\; l=0,\ldots,n+1,
$ 
and there exists
$\overline{\rho}\in(0,\frac{\rho_0}{2})$ such that
\begin{equation}\label{5.1}
\!\!\!\!\!\!\!\displaystyle \parallel u_{in}^\eta-u_{in} \parallel_{L^2 (-1,1)}<\frac{\sqrt{2}}{16Me^{\nu}}\eta,
\end{equation}
where $\nu$ is the nonnegative constant of \eqref{fsigni},
\begin{equation}\label{5.1bis}
M
:=\max_{x\in\overline{\Omega}_{\overline{\rho}}}
\Big\{\frac{\overline{u} (x)}{u_{in} (x)}\Big\}+1> \sup_{x\in \Omega_0
}
\Big\{\frac{\overline{u}_\eta (x)}{u_{in}^\eta (x)}\Big\},
\;
(\footnote{ We note that 
the auxiliary function $\Phi:[-1,1]\longrightarrow\R,$ defined in the following way
$$\Phi(x)=\begin{cases}
\dfrac{\overline{u}_{\eta} (x)}{u^{\eta}_{in} (x)},\;\; \text{ if }\,x\in \left(-1,1\right)\backslash \displaystyle\bigcup_{l=1}^n\left\{x_l\right\} \\
|\overline{u}^\prime_{\eta} (x)|, \;\:\,\text{ if }\,x=x_l,\, l=0,\ldots,n+1,
\end{cases}$$
is continuous on $[-1,1],$ since 
we have 
$\displaystyle \lim_{x\rightarrow x_l}\frac{\overline{u}_{\eta}(x)}{u^{\eta}_{in} (x)}=|\overline{u}_{\eta}^\prime (x_l)|,$\;\; for every $l=0,\ldots,n+1.$ 
})
\;\text{ and }\;\;
\parallel Mu^\eta_{in}\parallel_{L^2 \left(
(-1,1)\backslash 
\Omega
_{\overline{\rho}}
\right)}<\frac{\eta}{4}\,.
\end{equation} 
\end{itemize}
{\bf Step 1:} {\it Steering the system from $u_{in}+r_{in}
$ to $M u^\eta_{in}.$}
\,In this step, we represent the solution $u(x,t)$ 
of \eqref{u pres proof} as the sum of two functions $ w (x,t) $ and $h (x,t)$, 
which solve the problems in \eqref{w} in $\displaystyle(-1,1)\times(0,t_1),\, t_1>0,$  with
 the modified initial states:
 $$w|_{t = 0} \; =u_{in}^\eta\qquad \text{ and } \qquad h|_{t = 0} \; = r_{in}+\left(u_{in} -u_{in}^\eta\right).$$
Let us choose 
\begin{equation*}
\alpha (x,t) = 	\alpha_1 := \frac{\log M}{t_1}, \;\; (x,t) \in (-1,1)\times(0, t_1),
\end{equation*} 
for some $ t_1>0$. Applying the constant bilinear control  $ \alpha(x,t) = \alpha_1 > 0, \forall x\in (-1,1),$ on  the interval $(0, t_1),$ 
the solution of the first problem in \eqref{w} is given by
\begin{align}\label{w conv}
w (x,t_1)
= \; \sum_{p = 1}^\infty   &\left[\int_0^{t_1} e^{(\alpha_1 -\lambda_p^2 ) (t_1 -t)} \left(\int_{-1}^1 f (r,t,w (r, t))  \omega_p(r)
 \, dr \right) dt \right]\omega_p(x)
\nonumber\\
& +M
 \sum_{p = 1}^\infty   e^{-\lambda_p^2 t_1} \left(\int_{-1}^1 u^\eta_{in} (r) \omega_p(r)
  \, dr \right) \omega_p(x)
 \,,
\end{align}
where $\{-\lambda_p\}_{p\in\N}$
are the eigenvalues of the operator $(A_0,D(A_0)),$
and $\{\omega_p\}_{p\in\N}$ are the corresponding eigenfunctions, 
 that form a complete orthonormal system in $L^2(-1,1)$ (see also Proposition \ref{spectrum}).
By the strong continuity of the semigroup, see Proposition \ref{str cont}, 
we have that
 \begin{equation}\label{R}
 \sum_{p = 1}^\infty   e^{-\lambda_p^2 t_1} \left(\int_0^1 u^\eta_{in} (r) \omega_p(x)
  \, dr \right) \omega_p(x)\longrightarrow u^\eta_{in}(\cdot) \quad \text{ in } L^2(-1,1) \quad  \text{ as }\;\; t_1\rightarrow0.
 \end{equation}
Moreover, using
Proposition \ref{f in L2} 
 and H\"older's inequality we deduce
\begin{align}\label{H}
&\left\|\sum_{p = 1}^\infty   
\left[\int_0^{t_1} e^{(\alpha_1 -\lambda_p^2 ) (t_1 -t)}\left( \int_{-1}^1  
f (r,t,w (r, t))  
\omega_p(r)dr \right)
 dt \right]\omega_p(x)
 \right\|^2_{L^2 (-1,1)}
\nonumber\\ 
&=  \,
 \sum_{p = 1}^\infty \left| \int_0^{t_1} e^{(\alpha_1 -\lambda_p^2 ) (t_1 -t)}\left( \int_{-1}^1  
f (r,t,w (r, t))  \omega_p(r) \, 
dr \right)
 dt
\right|^2
\nonumber\\
&\leq  \,
 \sum_{p = 1}^\infty \left( \int_0^{t_1}  e^{2(\alpha_1 -\lambda_p^2 ) (t_1 -t)}dt
\right)\int_0^{t_1}\left|\int_{-1}^1f (r,t,w (r, t))  \omega_p(r) \, dr\right|^2\!dt 
\nonumber\\
 &\leq
 \sum_{p = 1}^\infty 
 e^{2\alpha_1t_1}t_1
\int_0^{t_1}\left|\int_{-1}^1f (r,t,w (r, t))  \omega_p(r) \, dr\right|^2\!dt 
=
M^2 
t_1 
\!\!\!\int_0^{t_1}\!\!\sum_{p = 1}^\infty\left|\int_{-1}^1f (r,t,w (r, t))  \omega_p(r)dr\right|^2\!\!\!\!dt\nonumber\\
&=M^2 
t_1 
\int_0^{t_1}
\int_{-1}^1f ^2(r,t,w (r, t)) drdt 
\leq 
 CM^2 e^{2k \vt t_1}
t_1\|u^\eta_{in}\|^{2\vt}_{1,a
}.
\end{align}
Making use of \eqref{w conv}-\eqref{H}, we have, 
as $t_1\rightarrow0$,
\begin{equation}\label{w final}
w (\cdot,t_1)\longrightarrow Mu^\eta_{in}(\cdot) \quad \text{ in } L^2(-1,1).
\end{equation}
{\it Now, we evaluate $ \|h (\cdot, t_1) \|_{L^2 (-1,1)}$. }
Let $t_1>0,$ multiplying by $ h$ both members of
 the equation in the second problem of (\ref{w}) and integrating by parts, since $\left[a(x)h_x(x,t)h(x,t)\right]_{-1}^{1}\leq0$ 
 we obtain
\begin{align*}
\int_{-1}^1 h^2 (x,t) dx 
&\leq\int_{-1}^1 
\left(r_{in}+\left(u_{in} -u_{in}^\eta\right)\right)^2 dx
+2\,\alpha_1 \int_0^{t_1}  \int_{-1}^1  h^2 dx dt\nonumber\\&
+2 \int_0^{t_1}\!\!\!\!\int_{-1}^1 \!\!\!(f(x,t,w+h)\!\! - \!\!f(x,t,w))h \,dx dt \nonumber\\&\leq 
 \int_{-1}^1 \left(r_{in}+\left(u_{in} -u_{in}^\eta\right)\right)^2 dx
+ 2(\alpha_1+\nu) \int_{0}^{t_1}  \int_{-1}^1  h^2 dx dt,\,\;\;\;t\in(0,t_1).
\end{align*}
 Thus, by Gronwall's inequality we have 
 $$\|h(\cdot,t)\|^2_{L^2(-1,1)}
\leq e^{2(\alpha_1+\nu){t_1}}
\|r_{in}+\left(u_{in} -u_{in}^\eta\right)\|_{L^2({-1},1)}^2,\;\;t\in (0,t_1),$$ then, for every $t_1\in(0,1),$
 \begin{equation}\label{h energy}
\|h(\cdot,t_1)\|_{L^2({-1},1)}\leq\|h\|_{C([0,t_1];L^2({-1},1))}
\leq 
Me^{\nu}\left(\|r_{in}\|_{L^2({-1},1)}+\|u_{in} -u_{in}^\eta\|_{L^2({-1},1)}\right).
\end{equation}
\noindent 
Since $u=w+h,$ by \eqref{w final}, \eqref{h energy} and \eqref{5.1}, there exists a small time $t_1\in(0,1)
$ such that 
\begin{equation}\label{5.10}
\|u (\cdot, t_1) - M u^\eta_{in} (\cdot)\|_{L^2({-1},1)}\leq \frac{\sqrt{2}}{8}\eta + Me^{\nu}\|r_{in}\|_{L^2({-1},1)}. 
\end{equation}
{\bf Step 2:} {\it Steering the system from $M u^\eta_{in}+\overline{r}_{in}
$ to $\overline{u}.$} \;
In this step, the solution $ u(x,t) $ of \eqref{u pres proof} is still represented  as the sum of  $ w (x,t) $ and $h (x,t)$, which solve the problems in \eqref{w} in $(-1,1)\times(t_1,T),$ with
 the modified initial states $M u^\eta_{in} $ instead of $ u_{in} $ and $\overline{r}_{in}(\cdot):=u (\cdot, t_1) - M u^\eta_{in} (\cdot)$ instead of $r_{in}(\cdot).$
By \eqref{5.1} it follows that
\begin{equation}\label{ass step}
\exists\, \delta^*>0\,:
\; \;\delta^*\leq\frac{\overline{u}_\eta
(x)}{M u^\eta_{in} (x)}
<1, \;\;\qquad \forall x \in \Omega_{\overline{\rho}}\,.
\end{equation} 
Owing to \eqref{ass step},
  assumption 
\eqref{ass} of Lemma \ref{lem preserving} is satisfied on $\,\Omega_{\overline{\rho}}\,$ with initial state $M u^\eta_{in}$ and target state $\overline{u}_\eta$\,. Then, adapting the proof of
Lemma \ref{lem preserving} 
we can conclude that there exists $T>t_1$ (such that $T-t_1>0$ is \textit{small}) and 
a bilinear control $\overline{\alpha}(x,t)=\frac{\tilde{\alpha}(x)}{T-t_1},\;\forall (x,t)\in(-1,1)\times(t_1,T),$ where
 $\tilde{\alpha}
\in C^2([-1,1])$
is very close in $L^2(-1,1)-$norm to the function\\ 
$\;\alpha_\eta (x)=\left\{ \begin{array}{ll}
\log \left( \frac{\overline{u}_\eta
 (x)}
{M u^\eta_{in} (x) } \right) \; & \text{ if } x \in \Omega_{\overline{\rho}},\\
0 \; &  {\rm elsewhere} \; {\rm in} \; [-1, 1]\,, \\
\end{array}
\right. 
$
such that
\begin{equation}\label{5.12}
\|u (\cdot, T) - \overline{u}_{\overline{\rho}}(\cdot)
 \|_{L^2(-1,1)}\leq \frac{\eta}{4}+\sqrt{2}\|u (\cdot, t_1) - M u^\eta_{in} (\cdot)\|_{L^2(-1,1)},
\end{equation}
$\text{ where }  \qquad\overline{u}_{\overline{\rho}}(x)=\left\{ \begin{array}{ll}
\overline{u}_\eta(x)
 \; & \text{ if } x \in \Omega_{\overline{\rho}},\\
Mu_{in}^\eta(x) \; &  {\rm elsewhere} \; {\rm in} \; [-1, 1]\,. \\
\end{array}
\right. $\\
Since $\!\|\overline{u}_\eta
-\overline{u}
\|_{L^2 (-1,1)}<\frac{\eta}{4},$ from \eqref{5.12}, 
\eqref{5.10} and \eqref{5.1bis} we obtain
\begin{align*}
\|u(\cdot,T) -\overline{u}(\cdot)\|_{L^2 (-1,1)}&\leq \|u(\cdot,T) -\overline{u}_{\overline{\rho}}
(\cdot)\|_{L^2 (-1,1)}\!\!
+
 \|\overline{u}_{\overline{\rho}} -\overline{u}_\eta
\|_{L^2 (-1,1)}\!\!+
\!\|\overline{u}_\eta
-\overline{u}
\|_{L^2 (-1,1)}\\
&\leq  \frac{\eta}{2}
+\sqrt{2}\|u (\cdot, t_1) - M u^\eta_{in} (\cdot)\|_{L^2(-1,1)}+\parallel Mu^\eta_{in}\parallel_{L^2 \left(
(-1,1)\backslash 
\Omega
_{\overline{\rho}}
\right)}\\
&\leq\eta+\sqrt{2}Me^{\nu}\|r_{in}\|_{L^2({-1},1)},
\end{align*}
from which the conclusion 
follows.
 $\qquad\qquad\qquad\qquad\qquad\qquad\qquad\qquad\quad\diamond$

\bigskip

{\bf References}

\bibliographystyle{elsarticle-num}

\end{document}